\documentclass[12pt]{article}

\usepackage{graphicx}
\usepackage[all,poly]{xy}

\title{Bayes factors and the geometry of discrete hierarchical loglinear models.}
\author{G\'erard Letac\thanks{Universit\'e Paul Sabatier} , H\'el\`ene Massam\thanks{Department of Statistics, York
University, Toronto, M3J 1P3, Canada. This author gratefully acknowledges support from  NSERC Discovery Grant A8947.}}
\date{}
\usepackage{latexsym}
\begin{document}
\maketitle

\def\<{\langle}
\def\>{\rangle}


\def\tr{\;\mathrm{tr}\;}
\def\<{\langle}
\def\>{\rangle}
\def\P{{\bf P}}
\def\E{{\bf E}}

\def\n{{\boldmath n}}
\def\b{{\boldmath b}}
\def \m{{\bf m}}
\def\z{{\boldmath z}}

\def\<{\langle}
\def\>{\rangle}
\newenvironment{pff}{\hspace*{-\parindent}{\bf Proof:}}{\hfill $\Box$
\vspace*{0.2cm}}
\newenvironment{pff.}{\hspace*{-\parindent}{\bf Proof}}{\hfill $\Box$}
\def\tr{\textmd{trace}\,}
\def\grad{\textmd{grad}\,}
\def\be{\begin{equation}}
 \def\ee{\end{equation}}
\def\beq{\begin{eqnarray*}}
 \def\eeq{\end{eqnarray*}}
 \def\tl{\triangleleft}
\newtheorem{theorem}{Theorem}[section]
\newtheorem{definition}{Definition}[section]
\newtheorem{subsecdef}{Definition}[subsection]
\newtheorem{lemma}{Lemma}[section]
\newtheorem{prop}{Proposition}[section]
\newtheorem{remarks}{Remarks}[section]
\newtheorem{remark}{Remark}[section]
\newtheorem{example}{Example}[section]
\newtheorem{perty} {Property}[section]
\newtheorem{cor} {Corollary}[section]

\textheight =22.5 cm
\textwidth =15 cm
\voffset =-0.5 in
\hoffset =0 in
\headheight =0 cm

\newtheorem{thm}{Theorem}[section]
\begin{abstract}
A standard tool for model selection in a Bayesian framework is the Bayes factor which compares the marginal likelihood of the data under two given different models. In this paper, we consider the class of hierarchical loglinear models for discrete data given under the form of a contingency table with multinomial sampling. We assume that  the Diaconis-Ylvisaker conjugate prior is the prior distribution on the loglinear parameters and the uniform is the prior distribution on the space of models. Under these conditions, the Bayes factor between two models  is a function of their prior and posterior normalizing constants. These constants  are functions of the hyperparameters $(m,\alpha)$ which can be interpreted respectively as marginal counts and the total count of a fictive contingency table.

We study the behaviour of the Bayes factor when $\alpha$ tends to zero. In this study two mathematical objects play a most important role. They are, first, the interior $C$ of the convex hull $\overline{C}$ of the support of the multinomial distribution for a given hierarchical loglinear model together with its faces and second, the characteristic function $J_C$ of this convex set $C$.
 We show that, when $\alpha$ tends to $0$, if the data lies on a face $F_i$ of $\overline{C_i},i=1,2$ of dimension $k_i$, the Bayes factor behaves like 
 $\alpha^{k_1-k_2}$. This implies in particular that when the data is in  $C_1$ and in $C_2$, i.e. when $k_i$ equals the dimension of model $J_i$, the sparser model is favored, thus confirming the idea of Bayesian regularization. 
 
 In order to identify the faces of $\overline{C}$, we need to know its facets. We give two new results. First, we identify a category of facets common to all hierarchical models for discrete variables, not necessarily binary. Second, we show that these facets are the only facets of $\overline{C}$ when the model is graphical with respect to a decomposable graph.

\vspace{.5cm}

\vspace{.5cm}

\end{abstract}

\newpage
\newpage
\section{Introduction}
We consider data given under the form of a contingency table representing the classification of $N$ individuals according to a finite set of criteria. We assume that the  cell counts in the contingency table follow a multinomial distribution. We also assume that the cell probabilities are modeled according to a hierarchical loglinear model (henceforth called hierarchical model).
The multinomial distribution for the hierarchical  model is a natural exponential family of the general form
$L(\theta)^{-1}\exp \<\theta,t\> \mu(dt)$
where $\mu$ is the generating measure and $L$ is its Laplace transform. The Diaconis-Ylvisaker \cite{dy79} conjugate prior has the general form 
\begin{equation}\label{FE}I(m,\alpha)^{-1}L(\theta)^{-\alpha}\exp (\alpha\<\theta,m\>) d \theta\end{equation}
 where $m$ and $\alpha$ are hyperparameters and $I(m,\alpha)$ is the normalization constant. 
Massam et al. \cite{mld09} have identified and studied the Diaconis-Ylvisaker  conjugate prior for the so called baseline constrained loglinear parametrization of the multinomial for hierarchical models.
This prior is a generalization of the hyper Dirichlet defined by Dawid and Lauritzen \cite{dl93} for graphical  models Markov with respect to decomposable graphs. Since decomposable graphical  models, and more generally graphical  models, form a subclass of the class of hierarchical models we will call this prior the generalized hyper Dirichlet. The hyper Dirichlet distribution is also used in discrete Bayesian networks as a prior for the cell parameters of the multinomial distribution of counts for the directed subgraphs formed by each discrete variable and its parents (see \cite{hgc95}). For the generalized hyper Dirichlet or the hyper Dirichlet, $\alpha$ is a positive scalar while  $m$ is a vector. The scalar $\alpha$ can be interpreted as the total sample size of a fictive contingency table and  $m$ can be interpreted as the vector of various marginal counts of the same  table. It is therefore traditional to take $\alpha$ small relatively to the total data count $N$. 
In this paper, we will use the loglinear parametrization  for the hierarchical model and  the generalized hyper Dirichlet as the prior, as defined in \cite{mld09}. 
\vspace{3mm}

In a Bayesian framework, the Bayes factor is one of the main tool for  model selection in the class of hierarchical models.
The aim of this paper is to study the behaviour of the Bayes factor for the  comparison of two hierarchical models $J_1$ and $J_2$ when $\alpha$ is very small, i.e., when $\alpha \rightarrow 0$. The motivation for this study is two-fold. First, it has been observed that as $\alpha \rightarrow 0$, in general, the Bayes factor will select the sparser model, that is the model with the parameter space of smallest dimension or equivalently the model with the least number of interactions. This is commonly  called the phenomenon of regularization, Second, Steck and Jaakkola \cite{sj02}, Proposition 1, have shown that, however, this is not always the case  and that, in fact, the behaviour of the Bayes factor between two Bayesian networks differing by one edge only depends upon a quantity which they call $d_{EDF}$, effective degrees of freedom, and which depends solely on the data. Comparing two such Bayesian networks is equivalent to comparing two graphical models on three variables, the saturated model and the model Markov with respect to the graph $A_3$, i.e., the two-link chain, with one conditional independence. It is therefore natural to seek a generalization  of the results in \cite{sj02} when  two arbitrary hierarchical models are considered.

Our aim is  to formally explain  when the sparser model is selected, when it is not and why. We also want to  develop tools to predict what the behaviour of the Bayes factor will be, for two given models. 

Since in the case of the Diaconis-Ylvisaker conjugate prior, the posterior probability of model $J$ given the data is equal to the ratio of the posterior and prior normalizing constants, we will be led to study the asymptotic behaviour, as $\alpha\rightarrow 0$, of the normalizing constant $I(m,\alpha)$  in (\ref{FE}). 
In this study, two important mathematical objects will surface. The multinomial distribution for a given hierarchical model $J$ is a natural exponential family. We denote by $C$ the interior of the convex hull $\overline{C}$ of the support of the measure generating this multinomial distribution. The convex polytope $\overline{C}$, together with its faces, is the first important mathematical object. The position of the data with respect to $C$, that is whether the data is in $C$ or on one of the faces of $\overline{C}$, will determine the behaviour of the Bayes factor. The second important mathematical object is the characteristic function $J_C$ of this polytope $C$; $J_C(m)$ is defined in the literature as the volume of the polar set of $C-m$ (see \cite{barv}). It is through $J_C$ that we will be able to find the asymptotic behaviour of $I(m,\alpha)$. Our central statistical result is that, as $\alpha \rightarrow 0$, the Bayes factor $B_{1,2}$ between two hierarchical models $J_1$ and $J_2$ behaves as follows:
 \begin{equation}
 \label{eq:bigBintro}
 B_{1,2}\sim D\alpha^{k_1-k_2}
 \end{equation}
where $D$ is a positive constant and $k_i,\;i=1,2$ are, respectively, the dimension of the face of $\overline{C_i}$ containing the data in its relative interior. When the data is in  both the open convex sets $C_i,\;i=1,2$, we have of course that
 $$B_{1,2}\sim D\alpha^{|J_1|-|J_2|}$$
 and this explains that in general the Bayes factor favours the sparser model since, in general for low-dimensional tables, the data is in the open polytope $C_i$. However with modern genetic or sociological data, we often deal with very sparse high-dimensional tables. In that case, the data may well be on a face of dimension  $k_i<|J_i|$. Then, as shown in \cite{sj02} for  three-factor models, the sparser model is not necessarily favoured by the Bayes factor. We do not consider, in this paper, the case $\alpha\rightarrow +\infty$ since in that case,  the behaviour of $I(m,\alpha)$ is well-known (see for example \cite{gideon} or \cite{haught88}).

\vspace{3mm} 
 
Here is a detailed description of the content of this paper. As mentioned above, we assume a multinomial  distribution for the counts and the generalized hyper Dirichlet prior for the loglinear parameters, as given in \cite{mld09}. However,  in order to  efficiently describe the geometry of $C$, we simplify the notation in \cite{mld09}. This new notation is given in Section 2.1. In Section 2.2, we derive a  characterization of the hierarchical  model which is close to that given in Proposition 3.1 of Darroch and Speed \cite{ds83}. Though not central to the derivation of our main result, this characterization  strengthens our understanding of the hierarchical model. In Section 2.3, we summarize \S 2 in \cite{mld09} and give a precise description of the measure $\mu$ generating the multinomial distribution, of its support and of the interior $C$ of the convex hull  $\overline{C}$ of this support.
 In Section 3, we describe properties of $J_C(m)$ and  $I(m,\alpha)$. Theorems
 \ref{jmlambda1}, \ref{imalpha} and \ref{jmlambda2} give, respectively, the expression of $J_C(m)$ in terms of the affine forms defining the facets of $\overline{C}$, the behaviour of $I(m,\alpha)$ when $m\in C$ and $\alpha\rightarrow 0$, and the behaviour of $J_C(z)$ as $z$ tends to the boundary of $C$ along a straight line. These results are used in Section 4 where we give our main statistical result, Theorem \ref{minimalpha} which yields equation (\ref{eq:bigBintro}). We thus have a precise description of the behaviour of the Bayes factor depending on the position of the data on the convex polytope $\overline{C}_i,\ i=\;1,2$ of the two models being compared. We show in Section \ref{steck} that our results comprise the results in \cite{sj02} as a particular case. In fact, we give a generalization of the concept of effective degrees of freedom to allow for the comparison of two arbitrary decomposable models. In Proposition \ref{generalEDF} using the generalized effective degrees of freedom, we give a quick and easy way to predict the behaviour of the Bayes factor.
 Since faces of $\overline{C}$ can only be obtained through the facets of $\overline{C}$, in Section 5, we return to the geometry of $\overline{C}$ and its facets. This set of facets has already been studied in the literature for certain binary hierarchical models (e.g. \cite{dl95} or \cite{hsul02}). In  Theorem \ref{generalfacets}, we describe a  category of facets common to all hierarchical models. 
 This constitutes a new result. For example, for the hierarchical model with four vertices $\{a,b,c,d\}$ and all three-way interactions $(abc),(bcd),(cda),(dab)$, no facets of $\overline{C}$ were known. In Corollary \ref{facedecomp}, we show that the special category of facets given in Theorem \ref{generalfacets} actually gives all the facets of $\overline{C}$ when the model is decomposable. We conjecture that this characterizes decomposable graphical models. Finally in Section 5.3, for the convenience of the reader, we present  some known results about the facets of $\overline{C}$ for graphical models Markov with respect to a cycle, using the notations of the present paper.

\section{Preliminaries}
\subsection{The notation}
While we keep the traditional notation as given in \cite{dl93} for cells and cell counts of the contingency table, we simplify the notation introduced in \cite{mld09} for the set of nonzero loglinear parameters. 

Let $V$ be a finite set of indices representing $|V|$ criteria. We assume that the criterion labelled by $v\in V$ can take values in a finite set $I_v$. We consider $N$ individuals classified according to these $|V|$ criteria. The resulting counts are gathered in a contingency table such that
$$I=\prod_{v\in V}I_{v}$$
is the  set of cells $i=(i_v,\;v\in V)$. If $D\subset V$  and $i\in I$ we write 
$i_D=(i_v, {v\in D})$ for the $D$-marginal cell. We write ${R}^I$ for the space of real functions $i\mapsto x(i)$ defined on $I.$ The element $x\in {R}^I$ is seen sometimes as a  vector, sometimes  as the function $i\mapsto x(i)$ on $I.$

Let $\mathcal{D}$ be a family of non empty subsets of $V$ such that $D\in \mathcal{D}$, $D_1\subset D$ and $D_1\not = \emptyset$  implies $D_1\in \mathcal{D}.$ In order to avoid trivialities we assume $\cup_{D\in \mathcal{D}}D=V.$  In the literature such a family $\mathcal{D}$ is called  a hypergraph (see \cite{laur96}) or an abstract simplicial complex (see \cite{efrs06}) or more simply the generating class (see \cite{eh85}). Following the notation introduced in \cite{ds83}, we denote by $\Omega_{\mathcal{D}}$ the linear subspace
of  $x\in {R}^I$  such that there exist functions $ \lambda_D\in {R}^I$ for $D\in \mathcal{D}$ depending only on $i_D$ and such that $x=\sum_{D\in \mathcal{D}}\lambda_D$, that is
\begin{eqnarray}
\label{omegad}\nonumber
\Omega_{\cal D}=\{x\in {R}^I:\;\exists\lambda_D\in {R}^I, D\in {\cal D}\;\mbox{such that}\; \lambda_D(i)=\lambda_D(i_D)\;\mbox{and}\; x=\sum_{D\in {\cal D}}\lambda_D\}
\end{eqnarray} 
The hierarchical  model generated by ${\cal D}$ is the set of probabilities $p=(p(i))_{i\in I}$ on $I$ such that $p(i)>0$ for all $i$ and such that $\log p\in \Omega_{\cal D}.$ It is convenient to write for $p$ in $\Omega_{\mathcal{D}}$
\begin{eqnarray}\label{H2}
\log p(i)=\lambda_{\emptyset}+\sum_{D\in \mathcal{D}} \lambda_D(i)
\end{eqnarray} 
where  $\lambda_{\emptyset}$ does not depend on $i$ and is thus a constant.  Needless to say the representation (\ref{H2}) is not  unique. 

We now introduce the notions we will need later to express the baseline constrained loglinear parameters used in the present paper. We first select a special element in  each 
$I_v$. For convenience we denote it $0$. By abuse of notation, we also denote $0$ in $I$ the cell with all its components equal to $0.$ This special element in $I_v$ is denoted $r_v$ in \cite{ds83} and $i^*$ in \cite{mld09}, but we find  the notation $0$ more convenient.  Actually the choice of the special element $0$ in each $I_v$ is arbitrary and does not affect our results. If $i\in I$ the support of $i$ is the subset of $V$ defined as $$S(i)=\{v\in V\ ; \ i_{v}\neq 0\}.$$  
We write
 \be \label{j} J=\{ j\in I,\ \ S(j)\in \mathcal{D}\}\ee 
 and note that since $\mathcal{D}$ does not contain the empty set, $J$ does not contain $0\in I$.
 This set $J\subset I$ is essential here and de facto defines the hierarchical  model. We introduce the important notation 
 $$j\tl i$$
  for $i\in I$ and $j\in J$ to mean that $S(j)$ is contained in $S(i)$ and that $j_{S(j)}=i_{S(j)}.$ Note that if $j,j'\in J$ and $i\in I$ we have 
 \begin{equation}\label{TL}j\tl j'\ \ \mbox{and}\ \ j'\tl i\Rightarrow j\tl i.\end{equation}
  Thus $\tl$ is in particular a partial ordering for $J$ but  we will never use the notation $i\tl i'$ for $i$ or $i'$ in  $I\setminus J$. Let us illustrate the notation above with an example. Let $V=\{a,b,c\}$, $\mathcal{D}=\{a,b,c,ab,bc\}$ and  $I_a=\{0,1,2\}=I_b$ and $I_c=\{0,1\}$. Thus $I$ has $3\times  3\times 2=18$ elements and 
 $$J=\{100,200,010,020,001,110,210,120,220,011,021\}$$ with 11 elements with respective supports $a,a,b,b,c,ab,ab,ab,ab,ac,ac.$ If $i=201$ the set of 
 $j$ in $ J$ such that $j\tl i$  is $\{200,001\}$ and if $i=211$ this set is $\{210,200,011,001,010\}.$
  
\subsection{The hierarchical model} 
 We now introduce  vectors which are fundamental for the description of the geometry of our problem. Let $(g_i)_{i\in I}$ and $(e_j)_{j\in J}$ be the canonical basis of ${R}^I$ and ${R}^J$ respectively. We endow these two spaces with their natural Euclidean structure. 
 For all $i\in I$ we define $f_i\in {R}^J$ by 
 \be\label{fi}
 f_i= \sum_{j\in J,\ j\tl i}e_j
 \ee 
 with in particular $f_0=0.$
 Let $H$ be the linear application $H:{R}^I\rightarrow {R}^J $  which sends the vector  $x=\sum_{i\in I}x(i)g_i$ of $R^I$ into
 $H(x)=\sum_{i\in I}x(i)f_i$ of $R^J.$
 Then 
 \be\label{H}H(x)=\sum_{j\in J}e_j\sum_{j\tl i}x(i).\ee The adjoint of $H$ is the linear application $H^*:{R}^J\rightarrow {R}^I$
  such that for all $\theta=\sum_{j\in J}\theta_je_j\in  {R}^J$  and all $x\in {R}^I$ one has
   $$\<H(x),\theta\>=\<\sum_{j\in J}e_j\sum_{j\tl i}x(i),\sum_{j\in J}\theta_je_j\>=\sum_{i\in I}x(i)\sum_{j\tl i}\theta_j=\<x,H^*(\theta)\>.$$
 As a consequence 
 \be\label{H*}
 H^*(\theta)=\sum _{i\in I}g_i\sum_{j\in J,\ j\tl i}\theta_j. 
 \ee 
 In other words the vector $H^*(\theta)$ of $R^I$ is the function $i\mapsto \sum_{j\in J,\ j\tl i}\theta_j.$ For instance for $j\in J$ one has 
 $$H^*(e_j)=\sum_{i\in I,\ j\tl i}g_i.$$
 Suppose that $x\in {R}^I$ is in the image of $H^*.$ The expression of $\theta\in {R}^J$  for a given $x=H^*(\theta)\in {R}^I$ is given in the following lemma.
 \begin{lemma}\label{exptheta}
 The mapping $H^*: R^J\mapsto R^I$ defined above is injective. If $x$ is in the image $\mathrm{im}(H^*)$ of $H^*$ we have $x=H^*(\theta)$ if and only if for all $j\in J$ 
 \be \label{IFx}\theta_j=\sum_{j'\in J\ ;\ j'\tl j}(-1)^{|S(j)|-|S(j')|}x(j').\ee 
In particular the vectors $(H^*(e_j))_{j\in J}$ are a basis of $\mathrm{im}(H^*).$

\end{lemma}

 \begin{pff}
 Let us first show the expression (\ref{IFx}) of $\theta_j,\  j\in J$. It follows from (\ref{H*}) that $x(j')=\sum_{j\tl j'}\theta_j$ and therefore (\ref{IFx}) is equivalent to
 $\theta_j=\sum_{j'\tl j}(-1)^{|S(j)|-|S(j')|}\sum_{j''\tl j'}\theta_{j''}$ which is equivalent to
 $$\theta_j=\sum_{j''\tl j}\theta_{j''}\sum_{j''\tl j'\tl j}(-1)^{|S(j)|-|S(j')|}.$$
 We therefore have to prove that
 for fixed $j''$ and $j$ in $J$ such that $j''\tl j'\tl j$,
  $$\sum_{j''\tl j'\tl j}(-1)^{|S(j)|-|S(j')|}=\left\{\begin{array}{ll} 1& \mbox{if}\; j=j''\\0&\mbox{if}\;j\neq j''.\end{array}\right .$$
 If $j=j''$ the result is trivially true. If $j\neq j''$ and $j''\tl j'\tl j$ with $j$ fixed, then $j'$ is entirely determined by its support $S(j)$ since $j'_{S(j')}=j_{S(j')}$. The principle of inclusion exclusion says that for $A\subset C$ where $C\neq A$ we have $\sum_{A\subset B\subset C}(-1)^{|B|}=0$. Applying this to $A=S(j'')$ and $C=S(j)$ gives the desired result and (\ref{IFx}) is proved. Note that (\ref{IFx}) implies that $H^*$ is injective and that the vectors $(H^*(e_j))_{j\in J}$ are independent.
 \end{pff}
 
 The following  proposition  characterizes ${\cal D}$ in terms of the linear application $H$ and its adjoint $H^*$. Its corollary describes the hierarchical model. 
 \begin{prop}\label{cod}
The space $\Omega_{\cal D}$ defined in (\ref{omegad}) is  the direct sum of the 1-dimensional space $K$ of constants and the image of $H^*$ in ${R}^J$
 $$\Omega_{\mathcal{D}}=K\oplus \mathrm{im} (H^*).$$
 The dimension of $\Omega_{\mathcal{D}}$ is 
$$d_{\cal D}=1+|J|\;\;\mbox{with}\;\;|J|=\sum_{D\in {\cal D}}\prod_{v\in D}(|I_v|-1).$$
\end{prop}

\begin{pff}
 We now show $K\oplus \mathrm{im} (H^*)=\Omega_{\mathcal{D}}.$ To see this, let us consider for $D\in \mathcal{D}$ the linear space $E_D$ of functions $i\mapsto \lambda_D(i)$ defined on $I$ and depending only on $i_D.$ This space is isomorphic to $R^{I_D}$ with the notation $I_D=\prod_{v\in D}I_v$ and therefore has dimension $|I_D|=\prod_{v\in D}|I_v|.$ 
 
 We now prove that a basis of the linear space $E_D$ is given by a vector generating the space $K$ of constants and by the set of the $|I_D|-1$ vectors of $R^I$ 
 $$\{H^*(e_j)\ j\in J,\ S(j)\subset D\}.$$
 
 To see this we observe that from  the definition of $H^*$  the value of the function $i\mapsto H^*(e_j)(i)$ is equal to $1$ if $j\tl i$ and to $0$  if not. If furthermore $S(j)\subset D$ this function $H^*(e_j)$ is an element of $E_D.$ This is checked by writing $i=(i_D,i_{D^c}):$ we have to show that $H^*(e_j)(i)$ does not depend on $i_{D^c}$. Consider the case  $H^*(e_j)(i_D,i_{D^c})=1$. This is saying  that $j\tl(i_D,i_{D^c})$ which implies that $S(j)\subset S(i)$ and $j_{S(j)}=i_{S(j)}.$ Recall our hypothesis  $S(j)\subset D.$ Now consider  $i'=(i_D,i'_{D^c}).$ Clearly $j\ \tl(i_D,i_{D^c})$ if and only if $j\ \tl(i_D,i'_{D^c})$, and this implies that $H^*(e_j)(i)=H^*(e_j)(i').$ Thus $H^*(e_j)\in E_D$ when $S(j)\subset D.$
 Now we recall that the vectors  $\{H^*(e_j)\ j\in J,\ S(j)\subset D\}$ are independent and that $K\cap \mathrm{im} (H^*)=\{0\}.$ Since the dimension of $E_D$ is  $(|I_D|-1) +1$  the claim is proved. 
 
 To complete the proof we use the fact that $\Omega_{\mathcal{D}}=\sum_{D\in \mathcal{D}}E_D$ (not a direct sum). Therefore
 $\Omega_{\mathcal{D}}= K\oplus \mathrm{im} (H^*).$  Note that $\mathrm{im} (H^*)$ can be seen as the subspace of the $x\in \Omega_{\mathcal{D}}$ such that $x(0)=0.$ The basis if $\mathrm{im}(H^*)$ being  the set of vectors $\{H^*(e_j)\ S(j)=D,\ D\in {\cal D}\}$ it then becomes clear that
 $$|J|=\sum_{D\in {\cal D}}\prod_{v\in D}(I_v-1).$$ The proposition is proved. \end{pff}

\begin{cor}\label{crux}The probability  $p=(p(i), i\in I)$ belongs  the hierarchical model generated by ${\mathcal{D}}$ if and only if there exist $\theta=\sum_{j\in J}\theta_je_j\in {R}^J$ and  a real number $\theta_0$ such that 
\be\label{logp}\log p=\theta_0+H^*(\theta)
\ee
 that is, for all $i\in I$, 
 $$\log p(i)=\theta_0+\sum _{j\in J\,\ j\tl i}\theta_{j}.$$ 
 Moreover,  $\theta\in {R}^J$ is uniquely defined by
 $$\theta_j=\sum_{j'\in J\ ;\ j'\tl j}(-1)^{|S(j)|-|S(j')|}\log p(j')$$
 and $\theta_0$ is uniquely defined by $e^{-\theta_0}=L(\theta)$, where
 \be \label{lap}
 e^{-\theta_0}=L(\theta)=1+\sum_{i\in I\setminus\{0\}}\exp(\sum_{j\in J\,\ j\tl i}\theta_{j}).
\ee
The discrete hierarchical  model generated by ${\cal D}$ is a manifold of dimension $|J|$. \end{cor}
The results of the corollary above are not new of course. The characterization  of the hierarchical model is close to that given in Proposition 3.1 of \cite{ds83} when one chooses what is called in that paper the substitution weight function for the averaging operator. 
\subsection{The multinomial distribution as a natural exponential family}  
  We will now use the proposition above to express the density of the multinomial distribution for the hierarchical  model generated by ${\cal D}$. We consider a contingency table with cells $i=(i_v,v\in V)\in I$ and cell counts $n=(n(i), \;i\in I)$  with $\sum_{i\in I}n(i)=N$ obtained from $N$ i.i.d. observations  of a multivariate Bernoulli variable with parameter $(p(i),\;i\in I)$, i.e. with distribution $\sum_{i\in {\cal I}}p(i)\delta_{g_i}.$ 
For $E\subset V$ we write $i_E\in I_E$ and $n(i_E)=\sum_{i'\in I; i_E=i'_E}n(i')$ for the $E$-marginal cell and $E$-marginal count respectively. For the particular case $E=S(j),\;j\in J$  we abbreviate $n(j_{S(j)})$ as
 \be \label{MC}t(j)=n(j_{S(j)}).\ee 
 Then, using (\ref{logp}), (\ref{H*}) and (\ref{H}) we have 
\begin{eqnarray*}
\sum_{i\in I} n(i)\log p(i)&=&\<\log p,n\>_{{R}^I}=N\theta_0+\<H^*(\theta),n\>_{{R}^I}=N\theta_0+\<\theta,H(n)\>_{{R}^J}\\&=&N\theta_0+\<\theta,t\>_{{R}^J}=N\theta_0+\sum_{j\in J} t(j)\theta_j\end{eqnarray*}
which, using (\ref{lap}) we rewrite
\be\label{rewrite}\prod_{i\in I}p(i)^{n(i)}=\frac{1}{L(\theta)^N}\exp\left(\sum_{j\in J}t(j)\theta_j\right)=\exp\left(\sum_{j\in J}t(j)\theta_j+N\theta_0\right).\ee 
The multinomial distribution for the model generated by ${\cal D}$ is therefore  a natural exponential family on ${R}^{J}$  and is generated by a discrete measure on ${R}^{J}$ whose Laplace transform is $L(\theta)^N.$   For $f_i$ as defined in (\ref{fi}), we have 
$L(\theta)=\sum_{i\in I}e^{\<\theta,f_i\>}$ and therefore $L$ is   the Laplace transform of the counting measure 
\be\label{cm}\mu=\sum_{i\in I}\delta_{f_i}\ee
 on the set of vectors $(f_i)_{i\in I}.$ This exponential family is concentrated on a bounded set of  ${R}^J$ and therefore the set of parameters $\theta$ for which $L$ is finite is the whole space ${R}^J.$ Hence the family is regular in the sense of Barndorff- Nielsen \cite{bn73} and Diaconis and Ylvisaker \cite{dy79}. Let  $C\subset{R}^J$ be the interior  of the convex hull of the set $(f_i)_{i\in I}$. In Corollary \ref{extremecor} below we show that the $(f_i)'s$ are the extreme points of its closure $\overline{C}.$ 
 \begin{prop}\label{extreme} Let $(e_j)_{j\in J}$ be the canonical basis of ${R}^J,$  let $(J_i)_{i\in I}$ be a family of subsets of $J$ such that $\cup_{i\in I}J_i=J$ and let $f_i=\sum_{j\in J_i}e_i,\; i\in I.$ 
  The extreme points of the convex hull $\overline{C}$ of the vectors $(f_i)_{i\in I}$ are the vectors $(f_i)_{i\in I}$ themselves.
 \end{prop}
 \begin{cor} \label{extremecor}The extreme points of the convex hull of the support of the measure $\mu$ as defined in (\ref{cm}) are the $f_i, i\in I$ as defined in (\ref{fi}).
 \end{cor}
 The corollary is obtained by taking $J_i=\{j\in J\ ;\ j\tl i\}$. Let us prove  Proposition \ref{extreme}.

 \begin{pff}
  Trivially any extreme point of $\overline{C}$ is an $f_i$ for some $i.$ Conversely 
let us show that for some given $i_0$, $f_{i_0}$ is  an extreme point of $\overline{C}$. Suppose  that there exist non negative numbers $(\lambda_i)_{i\in I}$ such that $\sum_{i\in I}\lambda_i=1$  and $f_{i_0}=\sum_{i\in I}\lambda_if_i.$ We are going to show that necessarily $\lambda_i=0$ if $i\neq i_0.$ By the definition of the $f_i$ 
\be\label{SL}f_{i_0}=\sum_{j\in J_{ i_0}} e_j=\sum_{i\in I}(\lambda_i\sum_{j\in J_i} e_j).\ee
We observe first that if $\lambda_i>0$ then $J_i\subset J_{ i_0}.$ If not there exists a $j_0\in J_i\setminus J_{i_0}$ and therefore  $\<f_{i_0},e_{j_0}\>=0.$ But (\ref{SL}) contradicts this since $\<f_{i_0},e_{j_0}\>\geq \lambda_i>0.$ 
Therefore,  writing
$$A(i_0)=\{i\in I; J_i\subset J_{i_0}\},$$ 
we must have $\lambda_i=0$ if $i\not \in A_{i_0}$ and $\sum_{i\in A(i_0)}\lambda_i=1.$ 
Then (\ref{SL}) becomes
\begin{eqnarray*}0 &=&(\sum_{i\in A(i_0)}\lambda_i)(\sum_{j\in J_{ i_0}} e_j)-\sum_{i\in A(i_0)}(\lambda_i\sum_{j\in J_ i} e_j)\\&=&
\sum_{j\in J_{ i_0}}e_j[(\sum_{i\in A(i_0)}\lambda_i)-(\sum_{i\in A(i_0),\ j \in J_i}\lambda_i)]\\&=&
\sum_{j\in J_{ i_0}}e_j(\sum_{i\in A(i_0), \ j\not \in J_ i}\lambda_i)\end{eqnarray*}
Since the $e_j$'s are independent, it follows that $\sum_{i\in A(i_0), \ j\not \in J_ i}\lambda_i=0$ and since the $\lambda$'s are nonnegative, this will imply $\lambda_i=0$ for each $i\neq i_0, i\in A(i_0)$ if we can show that there exists a $j_0$ such that $j_0\in J_{i_0}\setminus J_i$. This clearly true and therefore we conclude that $\lambda_i=0$ if $i\neq i_0.$ This proves that $f_{i_0}$ is an extreme point.
 \end{pff}
\subsection{The DY conjugate prior for the loglinear parameters}
As seen in Corollary \ref{crux}, the hierarchical  model generated by ${\cal D}$ can also be characterized by the set $J$ as defined in (\ref{j}). Moreover, from (\ref{rewrite}), we see that the multinomial distribution of $n$ for the  model $J$ can be written in terms of the marginal counts $t_J=(t(j), j\in J)$. It is the natural exponential family   with density, with respect to $\mu^{\otimes N}$, equal to
\be f(t|\theta, J)=C(n)\frac{\exp \<t_J,\theta\>}{L(\theta)^N}\ee
where $\mu$ is as in (\ref{cm}) and $C(n)$ is constant with respect to $\theta$. Following \cite{dy79}, the DY conjugate prior for $\theta$ indexed by $\alpha>0$ and by $m_J\in C$  is defined as the probability on ${R}^{J}$ with density 
$$\pi (\theta|m_J, \alpha, J)=\frac{1}{I_J(m_J,\alpha)}\times \frac{e^{\alpha\<\theta,m_J\>}}{L(\theta)^{\alpha}}$$
where $I_J(m,\alpha)$ is the normalizing constant. 
This family of priors is conjugate  and the  posterior probability of $\theta$ given the data  $n=(n(i))_{i\in I}$ in the contingency table is 
$$\pi(\theta|\frac{\alpha m_J +t_J}{\alpha+N},\alpha+N, J).$$ 
Let ${\cal H}$ denote the set of all hierarchical  models on the given set of variables. If we assume that the prior distribution on ${\cal H}$ is discrete, then the posterior distribution of $J$ given the data is
$$
g(J|t)=C(n)\frac {I_J(\frac{\alpha m_J +t_J}{\alpha+N},\alpha+N)}{I_J(m_J,\alpha)}\Big/\sum_{L\in {\cal H}}\frac {I_L(\frac{\alpha m_L +t_L}{\alpha+N},\alpha+N)}{I_L(m_L,\alpha)}
$$ 
In classical Bayesian model selection, the most probable models are selected by means of Bayes factors. More precisely, models are compared two by two by means of the Bayes factor $B_{1,2}$ between model $J_1$ and model $J_2$. In our framework
\be\label{BF}
B_{1,2}=\frac{I_2(m_2,\alpha)}{I_1(m_1,\alpha)}\times \frac{I_1(\frac{\alpha m_1 +t_1}{\alpha+N},\alpha+N)}{I_2(\frac{\alpha m_2 +t_2}{\alpha+N},\alpha+N)}\ee
where, for the sake of simplicity, $m,t,I$ are indexed by $i=1,2$ rather than by $J_1,J_2$ and where
$m_1$ and $m_2$ have been chosen in $C_1$ and $C_2$ respectively.
The aim of the present paper is to find  the limit of $B_{1,2}$ when $\alpha\rightarrow 0.$ If we assume that $n(i)>0$ for all $i\in I$, then  $t_k/N$ is in the interior of $C_k$ and under these circumstances the second factor in the right-hand side of (\ref{BF}) has the  finite limit $I_1(\frac{t_1}{N},N)/ I_2(\frac{t_2}{N},N).$ 
 For the first factor in (\ref{BF}), we will show that $I(m,\alpha)\sim_{\alpha\rightarrow 0} J_C(m)\alpha^{-|J|}$ where $J_C(m)$ will be studied in the next section. Thus when $\alpha\rightarrow 0$ the Bayes factor is equivalent to $$\alpha^{|J_1|-|J_2|}\frac{ J_{C_2}(m_2)}{J_{C_1}(m_1)}\times \frac{I_1(\frac{t_1}{N},N)}{I_2(\frac{t_2}{N},N)}.$$
 If we do not assume that $n(i)>0$ for all $i\in I$, then  $t_k/N$ might be on the boundary of $C_k$ for at least one $k=1,2$ and we will have to further study the behaviour of $I(m,\alpha)$ and $J_{C}(m)$. This is done
in the following section. 


\section{The limiting behavior of the prior normalizing constant} 
We give three fundamental theoretical results in this section. We assume that $m$ is in the interior of $C$, the convex hull of the measure $\mu$ as defined in (\ref{cm}).
Theorem \ref{jmlambda1} gives the general form of $J_C(m)$ in terms of the affine forms defining the facets of $C$. Theorem \ref{imalpha} gives the limit of $I(m,\alpha)$ when $\alpha\rightarrow 0$ and Theorem \ref{jmlambda2} describes the behaviour of $J_C((1-\lambda)y)$ when $y$ is on the boundary of $C$ and $\lambda \rightarrow 0$.
\subsection{The characteristic function of a convex set}
Given a finite dimensional real linear space $E$, let $E^*$ be its dual, that is, the space of all linear forms $\theta$ on $E$. We  write $\<\theta,x\>$ instead of $\theta(x)$ when $(\theta,x)\in E^*\times E.$  We fix a Lebesgue measure $d\theta$ on $E^*$  and a Lebesgue measure $dx$ on $E$ which must be compatible (this means that if $e$ is a basis of $E$ and $e^*$ is the corresponding dual basis of $E^*$ the product of the respective volumes of the two cubes built on $e$ and $e^*$ must be one). Needless to say when  $E={R}^n$ and $E^*=E$ and $\<.,.\>$ is the usual inner product and the Lebesgue measure is the usual one.  It will be however important in the sequel to distinguish between $E$ and $E^*$ and we therefore keep this notation. 

If $C\subset E$ is an open non empty convex set not containing a line, its polar set is  $$C^o=\{\theta\in E^*\ ;\ \<\theta,x\>\leq 1\ \forall x\in C\},$$ its support function $h_C:E^*\rightarrow (-\infty,\infty]$ is 
$$h_C(\theta)=\sup\{\<\theta,x\>\ ;\ x\in C\}$$ and its characteristic function is the function $m\mapsto J_C(m)$ defined on $ C$ by 
\be\label{JC}J_C(m)=\int_{E^*}e^{\<\theta,m\>-h_C(\theta)}d\theta.\ee 
We note that if $C$ contained a line, we would have $h_C(\theta)=\infty$ almost everywhere and $J_C\equiv 0.$ Faraut and Koranyi \cite{fk94}, p. 10  define $J_C $ when $C$ is an open convex salient cone.  In that case,  the polar set of $C$ is the convex cone
\be\label{PC}C^o=\{\theta\in E^*\ ;\ \<\theta,x\>\leq 0\ \forall x\in C\}\ee
and we have
$$h_C(\theta)=\left\{\begin{array}{ll}0&\mbox{if}\; \theta\in C^o\\+\infty& \mbox{if}\; \theta \not \in C^o\end{array}\right .$$
Let us mention here that when $C$ is a bounded set, $h_C(\theta)$ is finite for all $\theta\in E^*$. We also have the following important property of $J_C(.)$.
\begin{lemma}
\label{jfinite}
Let $C$ be an open convex set not containing a line and let $m\in C$. Then $J_C(m)$ is finite.
\end{lemma}

\begin{pff}
We  first give the proof for $m=0.$  Then, by assumption, $0$ is an interior point of $C$ and $h_C$ is always a strictly positive function. Assume that $E$ has a Euclidean structure and let $S(E)$ be its unit sphere. Recall that the function $\theta\mapsto h_C(\theta)$ is a continuous function. Thus, for $u=\theta/\|\theta\|\in S(E)$  we have the equality
$$h_C(\theta)/\|\theta\|=h_C(u)=\max\{\langle u,x\rangle\ ;\ x\in C\} .$$ 
Now the function $u\mapsto h_C(u)$ is continuous on the compact set $S(E)$: let $K>0$ be its minimum. The previous equality shows that
$$K\|\theta\|\leq h_C(\theta).$$ Thus if $n=\dim E$ we have 
$$\int _{E^*}e^{-h_C(\theta)}d\theta\leq \int _{E^*}e^{-K\|\theta\|}d\theta=C_n\int_0^{\infty}e^{-Kr}r^{n-1}dr<\infty$$
 where  $C_n=2\pi^{n/2}/\Gamma(n/2)$ is the area of $S(E).$
For the general case $m\neq 0$
we use the fact that  the support function of  $C-m$  satisfies
$h_{C-m}(\theta)=-\langle\theta,m\rangle +h_C(\theta).$   
\end{pff}

One can prove that $J_C(m)=\infty$ if $m\notin C.$ Another property of $J_C(m)$ is that when $C$ is an open convex set of $R^n$ not containing a line, the following formulas hold 
\begin{equation}
\label{eq:J}
J_C(m)=n!\mathrm{Vol}(C-m)^o=n!\int_{C^o}\frac{d \theta}{(1-\langle \theta, m\rangle)^{n+1}}
\end{equation}
For the first equality in (\ref{eq:J}), see \cite{barv} p. 207 and  \cite{hiriart} p. 243. For the second one, make the change of variable $\theta=\theta'/(1+\<\theta',m\>)$ in the integral $\int_{(C-m)^o}d\theta'.$ 

Computing $J_C(m)$ when $C$ is associated to an arbitrary hierarchical model is usually difficult except as we shall see in Section \ref{sec:decomposable}, when the model is a graphical decomposable model. Consider however  the following simple example:

\noindent {\it Example 1: the segment $(0,1)$.$\;$} Let $C=(0,1)\subseteq R$. In this case, $h_C(\theta)=\max (0,\theta)$ and for $0<m<1$ we have \be \label {JCB}J_C(m)= \int_{-\infty}^0e^{\theta m}d\theta+\int_0^{\infty}e^{\theta m-\theta}d\theta=\frac{1}{m}+\frac{1}{1-m}=\frac{1}{m(1-m)}.\ee 
Two more examples of $J_C(m)$ will be given after Theorem \ref{imalpha} below.

\vspace{3mm}

We now give a theorem that states that $J_C(m)$ is the ratio of polynomials where the denominator is equal to the product of the affine forms defining the facets of $\overline{C}$. This will be used in Section 5 to identify the facets of $\overline{C}$ for decomposable graphical models. We first need the following lemma which computes the characteristic function of a simplicial cone.
\begin{lemma}
\label{lemmavolume}
 Let $(x_1,\ldots,x_n)$ a basis of $E$ and let $(\xi_1,\ldots,\xi_n)$ be its dual basis in $E^*$ (that is $\<\xi_j,x_i\>=\delta_i^j$). Consider the simplicial cone $A$ of $E^*$ defined by
\begin{eqnarray*}
A&=&\{\theta=\theta_1\xi_1+\cdots+\theta_n\xi_n\; ; \; \theta_1> 0,\ldots,\theta_n> 0\}\\
&=&\{\theta\in E^*\; ;\;  \<\theta,x_1\>> 0,\ldots,\<\theta,x_n\>> 0\}
\end{eqnarray*}
and denote by $\mbox{Vol}(\xi_1,\ldots,\xi_n)$ the volume of the parallelotope 
$$\{\theta=\theta_1\xi_1+\cdots+\theta_n\xi_n\; ; \; 0\leq \theta_1\leq 1,\ldots,0\leq \theta_n\leq 1\}.$$
Then for all $x$ in  $-A^{o}\subset E$, that is the opposite of the dual cone of $A$, we have  $$\int_Ae^{-\<\theta,x\>}d\theta=\frac{\mathrm{Vol}(\xi_1,\ldots,\xi_n)}{\<\xi_1,x\>\ldots\<\xi_n,x\>}.$$
\end{lemma}
This lemma is elementary and is  obtained by writing $\theta$ in the $\xi$ basis and by making the change of variable from the coordinates of $\theta$ in the canonical basis of $R^n$ to the coordinates in the $\xi$ basis. 

 Recall that a \textit{facet} of a polytope $\overline{C}\subset R^n$ with a non empty interior is a face of dimension $n-1.$ More specifically a facet is the intersection of $\overline{C}$ with a supporting hyperplane of $\overline{C}$ which contains $n$ affinely independent points. 
\begin{theorem}\label{jmlambda1}
 Let $C\subset E$ be the non empty interior  of a bounded   polytope $\overline{C}$. Let  $m\in C.$ Then we have
 $$J_C(m)=\frac{N(m)}{D(m)}$$
  where
$D(m)=\prod_{k=1}^{K}g_k(m)$ is the product of affine forms $g_k(m)$ in $m$ such that $g_k(m)=0,\;k=1,\ldots,K$ define the facets of $\overline{C}$ and where $N(m)$ is a polynomial of degree $ <K.$ 
\end{theorem}

\begin{pff}
Let $\cal {E}$ be the set of extreme points of $\overline{C}$. By Corollary \ref{extremecor}, we know that $\overline{C}$ is the convex hull of $\cal{ E}$. Therefore for each $\theta$, there exists at least one $f\in {\cal E}$ such that $h_C(\theta)=\<\theta,f\>.$  Define the cone of influence of $f\in\cal {E}$ to be
 $$A(f)=\{\theta\in E^*;\ \langle \theta,x\rangle\leq \langle \theta,f\rangle\ \forall x\in \overline{C}\}=\{\theta\in E^*;\ h_C(\theta)= \langle \theta,f\rangle\}.$$ 
The cone of influence may be better visualized through its polar cone $A^o(f)$ which is contained in $E$ and is generated by $\overline{C}-f$. In other words, $f+A^o(f)$ is the support cone of $\overline{C}$ at its vertex $f$.

We now split $E^*$ into the  union of $A(f),\;f\in {\cal E}$ whose interiors are disjoint and intersections have measure zero. Indeed, for $f_i\in {\cal E},i=1,2$, 
$$A(f_1)\cap A(f_2)=\{\theta: \<\theta,f_1-g\>\geq 0\;\mbox{and}\;\<\theta,f_2-g\>\geq 0,\;\forall g\in \cal{E}\}$$
and therefore taking successively $g=f_1$ and $g=f_2$, we see that $A(f_1)\cap A(f_2)=\{\theta: \<\theta,f_1-f_2\>=0\}$ which is of measure $0$. Therefore, if we write
$$I_f(m)=\int_{A(f)}e^{\langle \theta, m\rangle -h_C(\theta)}d\theta=\int_{A(f)}e^{\langle \theta, m-f\rangle }d\theta$$ 
we have
$$
J_C(m)=\sum_{f\in \mathcal{E}}I_f(m)$$
In order to compute $I_f(m)$, we now split $A(f)$ into a union of closed simplicial cones $A_1(f),\ldots, A_{N_f}(f)$  with $n$ generators each,  with disjoint interiors  such that $\cup_{k=1}^{N_f}A_j(f)=A(f)$ and such that each generator of an $A_j(f)$ is a  generator of the cone $A(f).$
Each $A_j(f)$ is the intersection of $n$ half spaces given by
$$\{\theta\in E^*;\ \langle \theta,x^{(j)}_i\rangle\geq 0\},\;i=1,\ldots,n$$
for some vector $x_i^{(j)}$ of $E$ which is therefore an extreme generator of $A(f).$
Since the $A_j(f)$ are proper cones of $R ^n$,
$(x^{(j)}_1,\ldots,x^{(j)}_n)$ defines  a basis of $E$.  
The vector $f-m$ can be represented in this basis  as 
$$f-m=\sum_{i=1}^n(f_i^{(j)}-m_i^{(j)})x^{(j)}_i.$$
From Lemma \ref{lemmavolume}
\begin{equation}\label{ER}
I_f(m)=\sum_{j=1}^{N_f}\int_{A_j(f)}e^{-\langle \theta,f- m\rangle} d\theta=\sum_{j=1}^{N_f}\frac{\mbox{Vol}(\xi^{(j)}_1,\ldots,\xi^{(j)}_n)}{\prod_{i=1}^n(f_i^{(j)}-m_i^{(j)})}
\end{equation} 
where $(\xi^{(j)}_1,\ldots,\xi^{(j)}_n)$ is the dual basis of $(x^{(j)}_1,\ldots,x^{(j)}_n)$ and
$$f_i^{(j)}-m_i^{(j)}=\<\xi^{(j)}_i,f-m\>.$$
Reducing to the same denominator, we obtain that 
$$I_f(m)=\frac{P_f(m)}{\prod_{i=1}^M\langle \zeta_i,f-m\rangle}$$
where the $\zeta_i\in E^*$ are taken  among the $\xi^{(j)}_1,\ldots,\xi^{(j)}_n$ with $j=1,\ldots, N_f$ and where $P_f$ is a polynomial in $m$ with total degree $<M\leq nN_f.$  Note that for fixed $\zeta\in E^*$ the hyperplane of $E$ defined by 
$$H(f,\zeta)= \{m\in E;\ \langle \zeta,f-m\rangle=0\}$$ contains the extreme point $f.$ If $g$ is another extreme point of $\overline{C}$ and if $\langle \zeta_i,f\rangle=\langle \zeta_i,g\rangle$ then  $H(f,\zeta)=H(g,\zeta).$
This means that several factors of the denominator of $I_f(m)$ can also occur in $I_g(m)$ and therefore 
\begin{equation}
\label{jfactor}  
J_C(m)=\frac{N(m)}{\prod_{k=1}^Kg_k(m)}
\end{equation}
where $m\mapsto g_k(m),\;k=1,\ldots,K$ are distinct affine forms taken in the list of the $g_k(m)=\langle \zeta_i,f-m\rangle$ when $i$ and the extreme point $f$ vary, and where $N(m)$ is a polynomial in $m.$ Since, as a Laplace transform, $J_C$ is analytic in $C$ there is no point $m$ in  $C$ such that $g_k(m)=0.$ Therefore all facets  must be of the form $\overline{C}\cap \<\zeta_i,f-m\>=0$ for some $i\in \{1,\ldots,K\}$. Conversely every face $\overline{C}\cap \<\zeta_i,f-m\>=0$ is a facet: this is proved in the following general lemma.

\begin{lemma}\label{meu} Let $E$ be a $n$ dimensional space and let $\zeta_1,\ldots,\zeta_N$ generating  the dual space $E^*$, such that there exists $u\in E$ with $\<\zeta_j,u\>>0$ for all $j$ and such that $R_+\zeta_j$ is an extreme ray of the convex cone $A=\sum_{k=1}^NR_+\zeta_k$ for all $j=1,\ldots, N.$ If $B=A^o\subset E$ is the dual cone of $A$ then, for all $j=1,\ldots, N,$
 $$B\cap\{x\in E\ ;\ \<\zeta_j,x\>=0\}$$
is a facet of $B$. \end{lemma}
\begin{pff} Consider the hyperplane of $E^*$ defined by $P=\{\zeta\ ;\ \<\zeta,u\>=1\}.$
Without loss of generality we assume $\zeta_j\in P$ for all $j=1,\ldots, N.$ Thus $\zeta_1,\ldots,\zeta_N$ are the extreme points of the polytope $S=A\cap P.$ This polytope can be defined as the intersection of a finite number of $n-1$-dimensional half spaces $H_k\cap P, \;k=1,\ldots, K$ where $H_k=\{\zeta\in E^*\ ;\ \<\zeta,x_k\>\geq 0\}$ is a half-space in $E^*$ determined by $x_k\in E.$ Moreover, any particular extreme point $\zeta_j$ of $S$ is a face of $S$ of dimension $0$ and is therefore the intersection of $H_k\cap P, k\in I$ where $I\subseteq \{1,\ldots,K\}$ is of cardinality at least $n-1$. This is equivalent to saying that the linear system in $n-1$ unknown variables 
$$\{\zeta\in E^*\cap P\ ;\ \<\zeta,x_k\>= 0,\;k\in I\}$$ 
 has $\zeta_j$ as a unique solution and  we can therefore find $n-1$  vectors $(x_{k_j})_{j=1}^{n-1}$ which are independent. Since $B$ is the intersection of the half planes $\{x:\; \langle x, \zeta_j\rangle \geq 0\}, j=1,\ldots,N$ and since any $\zeta \in A$ can be written as a convex combination of $\zeta_j,j=1,\ldots,N$, the vectors $(x_{k_j})_{j=1}^{n-1}$ are in $B\cap\{x\in E\ ;\ \<\zeta_j,x\>=0\}$ and therefore define a facet of $B$. This completes the proof of Lemma \ref{meu}.
\end{pff}
The proof of Theorem \ref{jmlambda1} is also completed.
\end{pff}

\subsection{ The behaviour of $I(m,\alpha)$ as $\alpha\rightarrow 0$}
We have the following theorem.
\begin{theorem}\label{imalpha}
 Let $\mu$ be a positive measure on the finite dimensional linear space $E$ such that the interior $C$ of its closed convex support is not empty and is bounded. 
Denote by $L(\theta)=\int_Ee^{\<\theta,x\>}\mu(dx)$ its Laplace transform. For $m\in C$ and for $\alpha>0$ consider the Diaconis Ylvisaker integral
$$I(m,\alpha)=\int_{E^*} \frac{e^{\alpha \<\theta,m\>}}{L(\theta)^{\alpha}}d\theta.$$ Then 
\be \lim_{\alpha\rightarrow 0}\alpha^{n}I(m,\alpha)=J_C(m)\ee
where $n=\dim E.$ 
\end{theorem}
Let us note immediately that a remarkable feature of this result is that the limit $J_C(m)$ of $\alpha^{n}I(m,\alpha)$
depends on $\mu$ only through its convex support. For instance if $E=R$, the uniform measure on $(0,1)$ and the sum $\mu=\delta_0+\delta_1$ of two Dirac measures share the same $C=(0,1)$  and the same $J_C(m)=(m(1-m))^{-1}$.
We now need the following lemma.

\begin{lemma}\label{lemmaimalpha}
Let $\mu$ be  a bounded measure on some measurable space $\Omega$ and let $f$ be a positive, bounded  and measurable function on $\Omega$. Then we have
\begin{enumerate}
\item $||f||_p\rightarrow_{p\rightarrow \infty}||f||_{\infty}$
\item The function  $p\mapsto ||f||_p$ is either  decreasing on $(0,\infty)$ or there exists $p_0\geq 0$ such that it is  decreasing on $(0,p_0]$ and  increasing on $[p_0,+\infty)$.
\end{enumerate} 
\end{lemma} 

\begin{pff}
 Part 1 is well-known. Let us prove Part 2. Let $\nu$ be the image measure of $\mu$ by $f\mapsto \log f$, then
\begin{eqnarray*}
||f||_p&=&\Big(\int_{-\infty}^{+\infty} e^{px}\nu(dx)\Big)^{1/p}=\Big(\int_{-\infty}^{+\infty} e^{px}\nu_1(dx)\Big)^{1/p}C^{1/p}\\
&=&\exp \Big(\frac{1}{p}k_{\nu_1}(p)+\frac{1}{p}\log C\Big)
\end{eqnarray*}
where $\frac{\nu}{C}=\nu_1$ is a probability measure. Now we can easily verify through integration by parts that $k_{\nu_1}(p)=\int_0^p(p-t)k^{''}_{\nu_1}(t) dt$ and therefore
\begin{eqnarray*}
\frac{1}{p}k_{\nu_1}(p)&=&\int_0^p(1-\frac{t}{p})k^{''}_{\nu_1}(t) dt\\
\frac{d}{d p}(\frac{1}{p}k_{\nu_1}(p))&=&\frac{1}{p}\int_0^ptk^{''}_{\nu_1}(t) dt\\
\frac{d}{ dp}(\frac{1}{p}k_{\nu_1}(p)+\frac{1}{p}\log C)&=&\frac{1}{p^2}(\int_0^ptk^{''}_{\nu_1}(t) dt-\log C)
\end{eqnarray*}
Now since $k$ is strictly convex, the fonction $p\mapsto h(p)=\int_0^p tk^{''}_{\nu_1}(t) dt-\log C$ is continuous and increasing. Therefore either $h$ is negative for all $p$ or $h$ is negative until its unique zero $p_0$  and then it is positive. This proves the lemma. 
\end{pff}
 
\begin{pff} (\textsc{of Theorem \ref{imalpha}.}) In the integral $\alpha^nI(m,\alpha)$ we make the change of variable $y=\alpha\theta$ and we obtain 
$$\alpha^nI(m,\alpha)=\int_{E^*}\frac{e^{\langle y,m\rangle}}{L(y/\alpha)^{\alpha}}dy.$$
We now apply the last lemma to $\Omega=\overline{C},$ to the bounded measure $\mu,$  to the function $f(x)=e^{\langle y,x\rangle}$ for some fixed $y\in E^*$
and to $p=1/\alpha.$ Denote by $S$ the  support of $\mu$. One easily sees that the support function of $C$ satisfies 
$$h_C(\theta)=\sup\{\langle \theta,x\rangle\ ;\ x\in C\}=\max\{\langle \theta,x\rangle\ ;\ x\in S\}$$ since $C$ is the interior of the convex hull of $S.$ As a consequence the  essential sup of $f$ is $e^{h_C(y)}$ and we get 
$$\lim_{\alpha\rightarrow 0}L(y/\alpha)^{\alpha}=e^{h_C(y)}.$$ 
Furthermore, by Lemma \ref{lemmaimalpha}, the function $p\mapsto \|f\|_p$ is monotonic for $p$ big enough. If $p\mapsto ||f||_p$ is increasing, $\frac{1}{||f||_p||}$ is decreasing and then by the monotone convergence theorem
\begin{eqnarray*}
\lim_{\alpha\rightarrow0}\int_{E^*}\frac{e^{\langle y,m\rangle}}{L(\frac{y}{\alpha})^{\alpha}}dy&=& \int_{E^*}\frac{e^{\langle y,m\rangle}}{\lim_{\alpha\rightarrow 0}L(\frac{y}{\alpha})^{\alpha}}dy=
\int_{E^*}\frac{e^{\langle y,m\rangle}}{\lim_{p\rightarrow \infty}||e^{\langle y,m}||_p}dy\\&=&\int_{E^*}\frac{e^{\langle y,m\rangle}}{||e^{\langle y,m}||_{\infty}}dy\int_{E^*}e^{\langle y,m\rangle-h_C(y)}dy=J_C(m)\nonumber
\end{eqnarray*}
If $p\mapsto ||f||_p$ is  decreasing,  $p\mapsto 1/||f||_p$ is increasing. In order to show that we can invert the order of limit and integration and to apply the monotone convergence theorem as we did in the previous case, we need here insure that $\int_{E^*}e^{\langle y,m\rangle-h_C(y)}dy$ is finite: Lemma \ref{jfinite} shows that it is true. The proof is complete.
\end{pff}
We now give two more examples of functions $J_C(m)$ which we compute using Theorem \ref{imalpha}.

\noindent {\it Example 2: $\overline{C}$ is the simplex.$\;$}  Let $e_0=0$ and $(e_1,\ldots,e_n)$ be the canonical basis of $R^n$. Let $C$ be the interior of the convex set generated by $e_0,\ldots, e_n.$ Then
$C$ is the set of $m\in R^n$ such that $m=\sum_{j=0}^n\lambda_j e_j$ for some unique positive $\lambda_0,\ldots,\lambda_n$ satisfying $\lambda_1+\cdots+\lambda_n<1$.
In this case 
$$J_{C}(m)=\frac{1}{m_1m_2\ldots m_n(1-m_1-\cdots-m_n)}.$$
This result can be obtained by computing $I(m,\alpha)$ for $\mu(x)$ in (\ref{cm}) equal to $\mu=\delta_{e_0}+\sum_{i=1}^n \delta_{e_i}$. Using elementary methods of integration, we find that
$$I(m,\alpha)=\int_{R^n}\frac{e^{\alpha\<\theta,m\>}}{(1+\sum_{i=1}^ne^{\theta_i})^{\alpha}}=\int_{R^n}\frac{\prod_{i=1}^n e^{\alpha m_i \theta_i}}{(1+\sum_{i=1}^n e^{\theta_i})^{\alpha}}\prod_{i=1}^n d \theta_i=\frac{\prod_{i=0}^n \Gamma(\alpha m_i)}{\Gamma(\sum_{i=0}^n \alpha m_i)}$$
where $m_0=1-\sum_{i=1}^n m_i.$ Using  $z\Gamma(z)=\Gamma(1+z)\rightarrow_{z\rightarrow 0}= 1$ we immediately obtain that
$$J_C(m)=\lim_{\alpha\rightarrow 0}\alpha^n I(m,\alpha)=\frac{1}{\prod_{i=0}^n m_i}.$$
 
\noindent {\it Example 3: $C$ for the graphical model $\stackrel{a}{\bullet}-\stackrel{b}{\bullet} -\stackrel{c}{\bullet}$.$\;$} For simplicity, we will assume that the variables $a,b,c$ are binary so that $m=(m_j, j\in J)$ where $J$ is defined as in (\ref{TL}) can be written $m=(m_D, D\in {\cal D})$ where
${\cal D}=\{a,b,c,ab, bc\}.$ We shall generalize this example in Section 5. 
From formula (4.8) in \cite{mld09}, we know that
 \begin{eqnarray*}I(m,\alpha)&=&\Gamma(\alpha(1-m_a-m_b+m_{ab}))\Gamma(\alpha(m_a-m_{ab}))\Gamma(\alpha(m_b-m_{ab}))\\&&\times\Gamma(\alpha(m_{ab}))
 \Gamma(\alpha(1-m_b-m_c+m_{bc}))\Gamma(\alpha(m_b-m_{bc}))\Gamma(\alpha(m_c-m_{bc}))\Gamma(\alpha(m_{bc}))\\&&\times\frac{1}{\Gamma(\alpha m_b)\Gamma(\alpha(1-m_b))}\end{eqnarray*}
 and therefore using  $z\Gamma(z)=\Gamma(1+z)\rightarrow_{z\rightarrow 0} 1$  again we obtain that
 \begin{eqnarray*}
 \lim_{\alpha\rightarrow 0}\alpha^5I(m,\alpha)&=&J_C(m)=\frac{m_b(1-m_b)}{m_{ab}m_{bc}}\\&\times&\frac{1}{(1-m_a-m_b+m_{ab})(m_a-m_{ab})(m_b-m_{ab})}\\
&\times &\frac{1}{(1-m_b-m_c+m_{bc})(m_b-m_{bc})(m_c-m_{bc})}
\end{eqnarray*}

\subsection{ The behaviour of $J_C(\lambda m+(1-\lambda)y)$ when $y\in \overline{C}\setminus C$ and $\lambda\rightarrow 0$}
In practice the choice of the hyperparameters $m$ and $\alpha$ is ours and for a given model $J$, it is traditional to take $m=(m_j,\;j\in J)$ to be the vector of $J$-marginal counts in a  fictive contingency table with cell counts all equal and equal to $\frac{1}{|I|}$. In any case, as long as all fictive cell counts are positive, $m$ belongs to the open set $C$ and the behaviour of $I(m,\alpha)$ is given by Theorem \ref{imalpha}.

 When studying the Bayes factor, we will have to consider the case where the data belongs to the boundary $\overline{C}\setminus C=\partial C$ of $C$, that is to a face of $\overline{C}$.
  To do so, we will need to describe  the behaviour of  $J_C(z)$ as $z$ approaches the boundary of $C$ along a straight line. This is done in the following theorem. Without loss of generality, we assume that $m=0$ so that
 $J_C(\lambda m+(1-\lambda)y)$=$J_C((1-\lambda)y)$. 
\begin{theorem} 
\label{jmlambda2} 
 Let $C$ be a polytope $\subset E$ with $\dim E=n$ and such that $0$ is in the interior of $C.$ Let $y\in \partial C,$ let $F$ be the face of $C$ containing $y$ in its relative interior and let $k$ be the dimension of $F$. Then when $\lambda\rightarrow 0$
$$\lim_{\lambda\rightarrow 0}\lambda^{n-k}J_C((1-\lambda)y)=D,$$ 
where $D$ is a positive constant.
\end{theorem} 

\begin{pff}  
From (\ref{eq:J}) we have
\begin{equation}
\label{basic}
 \frac{J_C((1-\lambda) y)}{n!}=\int_{C^o} \frac{d\theta}{(1-(1-\lambda)\<\theta,y\>)^{n+1}}\;.
 \end{equation}
In order to study the behaviour of this last integral when $\epsilon\rightarrow 0$, we are  going to build a parametrization of $\overline{C^o}$ which gives a special role to $\widehat{F}$, the face of $\overline{C^o}$ dual to the face $F$ of $\overline{C}$ containing $y$ in its interior.

Let $\mathcal{E}$ the set of extreme points of $\overline{C}$ and $\mathcal{I}\subset \mathcal{E}$ the set of extreme points of $F.$ To $F$ we associate the dual face of $\overline{C^o}$ defined by
\begin{equation}
\label{FC}\widehat{F}=\{\theta\in\overline{ C^o}\ |\ \<\theta,f\> =1\ \forall f\in \mathcal{I}\}.
\end{equation} 
It is a classical result (see \cite{barv}) that $\widehat{F}$ has dimension $n-k-1.$
Let us now observe that we have an equivalent representation of $\widehat{F}$ in (\ref{FC}) as
\begin{equation}
\label{FA}
\widehat{F}=\{\theta\in C^o\ |\ \<\theta,y\> =1\}.\end{equation}
Indeed, since  $y$ is in the relative interior of $F$ we write
$$y=\sum_{f\in \mathcal{I}}\lambda_f f$$
where $\lambda_f>0$ and $\sum_{f\in \mathcal{I}}\lambda_f=1.$ Here $\lambda_f>0$ is important in the argument to follow. Clearly $\widehat{F}\subset \{\theta\in C^o\ ;\ \<\theta,y\> =1\}.$ Conversely if $\<\theta,y\> =1$
then $\sum_{f\in \mathcal{I}}\lambda_f(1-\<\theta,f\>)=0.$ If furthermore $\theta\in \overline{C^o}$ we have $1-\<\theta,f\>\geq 0$ and therefore $1-\<\theta,f\>=0$
which shows $\widehat{F}\supset \{\theta\in \overline{C^o}\ ;\ \<\theta,y\> =1\}$ and proves (\ref{FA}).

Next, for  $\epsilon>0$ small, we consider the following approximation $\widehat{F}_{\epsilon}$ of $\widehat{F}$ 
\begin{equation}
\label{FB}
\widehat{F}_{\epsilon}=\{\theta\in \overline{C^o}\ ;\ \<\theta,y\> =1-\epsilon\}.
\end{equation} which is a $n-1$ dimensional convex subset of $\overline{C^o}$  and we want to prove that $$\mathrm{vol}_{n-1}\widehat{F}_{\epsilon}\sim c\epsilon ^{k}$$ for some positive constant $c$. Using (\ref{FA}) and(\ref{FC}), we can rewrite (\ref{FB}) as
\begin{equation}\label{FD}\widehat{F}_{\epsilon}=\{\theta\in \overline{C^o}\ ;\ \sum_{f\in \mathcal{I}}\lambda_f(1-\<\theta,f\>) =\epsilon\ \}.\end{equation} 
To show $\mathrm{vol}_{n-1}\widehat{F}_{\epsilon}\sim c\epsilon ^{k}$ we parametrize $\widehat{F}_{\epsilon}$ as follows:   let $\theta\mapsto x=\varphi(\theta)$ be the affine map from $E^*$ to $R^{\mathcal{I}}$ defined by 
\begin{equation}\label{SC}
x_f=\lambda_f(1-\<\theta,f\>),\;\;f\in {\cal I}
\end{equation}
which is equivalent to $\<\theta,f\>=1-\frac{x_f}{\lambda_f}$. The set  $S_{\epsilon}=\varphi(\widehat{F}_{\epsilon})$  is therefore the intersection of the simplex 
\begin{equation}
\label{MS}
\{x\ inR^{\mathcal{I}}\ ;\ x_f\geq 0 \ \forall f\in \mathcal{I},\ \sum_{f\in \mathcal{I}}x_f=\epsilon\}
\end{equation} 
and of the convex set $\varphi(\overline{C^o})$ which is contained in the affine manifold $\varphi(E^*)\subset  R^{\mathcal{I}}.$ If $x\in S_{\epsilon}$ then its preimage by $\varphi$ is the set
$$\varphi^{-1}(x)=\{\theta\in E^*\ ;\ \<\theta,f\>=1-\frac{x_f}{\lambda_f}\ \forall f\in \mathcal{I}\}$$ which is an affine subspace of $E^*$ parallel to the linear space 
\begin{equation}
\label
{MH}H_0=\{\theta\in E^*\ ;\ \<\theta,f\> =0\ \forall f\in \mathcal{I}\}
\end{equation} 
which has dimension $n-k-1$ since $F$ has dimension $k$. 
 As a result we can write $\widehat{F}_{\epsilon}$ as the following union of disjoint sets
\begin{equation}
\label{paraFe}
\widehat{F}_{\epsilon}=\cup_{x\in S_{\epsilon}}(\varphi^{-1}(x)\cap C^o)
\end{equation}
which is saying that $\widehat{F}_{\epsilon}$ can be parametrized by $(x,z)$ where $x\in S_{\epsilon},$ a convex set of dimension $k,$ and where $z\in \varphi^{-1}(x)\cap \overline{C^o},$  a convex set of  dimension $n-k-1.$ The bijection $\theta\mapsto (x,z)$ is the restriction to $\widehat{F}_{\epsilon}$ of the affine map $\varphi$ and therefore its Jacobian $K$ such that $d\theta=Kdxdz$ is a constant. 
$$\mathrm{vol}_{n-1}\widehat{F}_{\epsilon}=\int_{\widehat{F}_{\epsilon}}d\theta=K\int_{S_{\epsilon}}\left(\int_{\varphi^{-1}(x)\cap \overline{C^o}}dz\right)dx$$
If we fix $\epsilon x^0$ in the simplex (\ref{MS}), then 
the behavior of $\int_{\varphi^{-1}(\epsilon x^0)\cap \overline{C^o}}dz$ is easy to describe since 
$\lim_{\epsilon } \varphi^{-1}(\epsilon x^0)\cap \overline{C^o}=\widehat{F}$ in the sense of polytopes, which implies 
$$\lim_{\epsilon \rightarrow 0}\int_{\varphi^{-1}(\epsilon x^0)\cap C^o}dz=\lim_{\epsilon \rightarrow 0 }\mathrm{vol}_{n-k-1}(\varphi^{-1}(\epsilon x^0)\cap \overline{C^o})=\mathrm{vol}_{n-k-1}(\widehat{F}).$$
Let us now observe that $0$ is an extreme point of $\varphi(\overline{C^o}).$ If not there exist $x=\varphi(\theta)$ and $x'=\varphi(\theta')$ with $\theta$ and $\theta'\in C^o$ such that $x+x'=0$, that is, for all $f\in {\cal I}$
$$1-\lambda_f\<\theta,f\>+1-\lambda_f\<\theta',f\>=2-\lambda_f[\<\theta,f\>+\<\theta',f\>]=0.$$
Since $0\leq \lambda_f\leq 1$, this in turn
 implies $\lambda_f=1$ and $\<\theta,f\>+\<\theta',f\>=2.$ Since $\<\theta,f\>$ and $\<\theta',f\>$ are $\leq 1$ this implies $x_f=x'_f=0$ for all $f\in\mathcal{I},$ a contradiction. Now we use the fact that $C^o$ is a polytope and so is $\varphi(\overline{C^o})$ which has dimension $k.$ For $\epsilon$ small enough (say $0<\epsilon\leq \epsilon_0$) the intersection $S_{\epsilon}$ of the simplex given in (\ref{MS}) with $\varphi(\overline{C^o})$ coincides with the intersection of the simplex with the support cone of $\varphi(\overline{C^o})$ at its vertex $0.$ Since a cone is invariant by dilations we can claim that there exists a number $c_1>0$ such that for $0<\epsilon\leq \epsilon_0$ we have 
$$\mathrm{vol}_{k}(S_{\epsilon})=c_1\epsilon ^k.$$ Finally 
\begin{equation}
\label{int}
\mathrm{vol}_{n-1}\widehat{F}_{\epsilon}\sim c_1\; K\; \mathrm{vol}_{n-k-1}(\widehat{F})\epsilon ^k.
\end{equation}
The parametrization of $\theta$ in (\ref{basic}) is therefore $(x,z, \epsilon)$ where $(x,z)$ is as given in (\ref{paraFe}) and 
the range of $\epsilon$ is such that, for that range, $F_{\epsilon}$ describes all of $\overline{C^o}$. We note that the bounded function $\mathrm{vol}_{n-1}\widehat{F}_{\epsilon}=f(\epsilon)$ is zero if $\epsilon$ is big enough  since then  $\widehat{F}_{\epsilon}$ becomes empty and, of course, $\mathrm{vol}_{n-1}\widehat{F}_{0}=\mathrm{vol}_{n-1}\widehat{F}$. Let $b$ be such that $f(\epsilon)=0$ when $\epsilon>b$. When $\epsilon$ varies from $0$ to $+\infty$, $\widehat{F_{\epsilon}}$ generates all of $\overline{C^o}$. Then, following (\ref{int}), equation (\ref{basic}) becomes
\begin{eqnarray*}
\int_{\overline{C^o}} \frac{d\theta}{(1-(1-\lambda)\<\theta,y\>)^{n+1}}&=&\int_0^\infty \frac{\mathrm{vol}_{n-1}\widehat{F}_{\epsilon}d\epsilon}{(1-(1-\lambda)(1-\epsilon))^{n+1}}\\
&=&\int_0^\infty \frac {f(\epsilon) d\epsilon}{(1-(1-\lambda)(1-\epsilon))^{n+1}}
\end{eqnarray*}

Using $f(\epsilon)\sim c\; \epsilon^k$ we will now show that 
\begin{equation}
\label{MF}
\lim_{\lambda\rightarrow 0}\lambda^{n-k}\int_0^\infty \frac{f(\epsilon)d\epsilon}{(1-(1-\lambda)(1-\epsilon))^{n+1}}=c\; B(k+1,n-k),
\end{equation} and this will conclude the proof. 
To derive (\ref{MF}), we first show that for $0<a<b$
\begin{enumerate}
\item[(1)] $\lambda^{n-k}\int_0^a\frac{\epsilon^k d\epsilon}{(\lambda+\epsilon-\lambda \epsilon)^{n+1}}\rightarrow_{\lambda\rightarrow 0} B(k+1,n-k)$
\item[(2)] $\lim_{\lambda\rightarrow 0}\lambda^{n-k}\int_a^b\frac{ d\epsilon}{(\lambda+\epsilon-\lambda \epsilon)^{n+1}}=0$
\end{enumerate}
Statement (1) is shown by the change of variable $\epsilon=\lambda t$ and the theorem of dominated convergence. Indeed, for $0<\lambda\leq \lambda_0<1$, we have 
\begin{eqnarray*}
\lambda^{n-k}\int_0^a\frac{\epsilon^k d\epsilon}{(\lambda+\epsilon-\lambda \epsilon)^{n+1}}&=&\int_0^{a/\lambda}\frac{t^k}{(1+t-\lambda t)^{n+1}}\le \int_0^{a/\lambda}\frac{t^k}{(1+t-\lambda_0 t)^{n+1}}
\end{eqnarray*}
which tends to $\frac{1}{(1-\lambda_0)^{k+1}}B(k+1,n-k)$ when $\lambda\rightarrow 0$. Since this is true for any $\lambda_0>0$, statement (1) follows.

Statement (2) is obvious since
$\int_a^b\frac{ d\epsilon}{(\lambda+\epsilon-\lambda \epsilon)^{n+1}}<\int_a^b\frac{ d\epsilon}{\epsilon^{n+1}}$ is finite.
\noindent Next, fix $\delta >0.$ There exists $a<b$ such that $|\frac{f(\epsilon)}{\epsilon^k}-c|\leq \delta$ if $0<\epsilon \leq a$. Writing this as $-\delta\epsilon^k<f(\epsilon)-c\epsilon^k<\delta\epsilon^k$, integrating and using (1) yields  
$$\limsup _{\lambda\rightarrow 0}\left|\frac{1}{B(k+1,n-k)}\int_0^a \frac{f(\epsilon)d\epsilon}{(1-(1-\lambda)(1-\epsilon))^{n+1}}-c\right|\leq \delta.$$
 Since $f$ is bounded (2) implies that $$\limsup _{\lambda\rightarrow 0}\lambda^{n-k}\int_a^{+\infty} \frac{f(\epsilon)d\epsilon}{(1-(1-\lambda)(1-\epsilon))^{n+1}}=\limsup _{\lambda\rightarrow 0}\lambda^{n-k}\int_a^b \frac{f(\epsilon)d\epsilon}{(1-(1-\lambda)(1-\epsilon))^{n+1}}=0$$
Thus for all $\delta>0$ we have $$\limsup _{\lambda\rightarrow 0}\left|\frac{1}{B(k+1,n-k)}\int_0^{\infty} \frac{f(\epsilon)d\epsilon}{(1-(1-\lambda)(1-\epsilon))^{n+1}}-c\right|\leq \delta,$$ which implies (\ref{MF}).
\end{pff}

\section{The limiting behaviour of the Bayes factor}
Let us recall that, under the uniform distribution on the class of hierarchical models, the Bayes factor between two models $J_1$ and $J_2$ is equal to
\begin{eqnarray*}
B_{1,2}&=&\frac{I_1(\frac{\alpha m_1+t_1}{\alpha+N},\alpha+N)I_2(m_2,\alpha)}{I_2(\frac{\alpha m_2+t_2}{\alpha+N},\alpha+N)I_1(m_1,\alpha)}
\;.
\end{eqnarray*}
where $t_i=t_{J_i}=(t(j),\;j\in J_i),\;i=1,2.$
\subsection{The  case where the data is in the interior $C$ of $\overline{C}$}
The data is of course given in the form of a contingency table with cell counts $n=(n(i),\;i\in I)$. In this subsection, we consider the case where the data, which appears under the form $t_i$ in models $J_i,$ belongs to $C_i,\;i=1,2$ so that $I_i(\frac{t_i}{N},N),\;i=1,2$ are finite and positive.
In this case, as $\alpha \rightarrow 0$, from Theorem \ref{imalpha}, we know that, as $\alpha\rightarrow 0$,
\begin{eqnarray}\label{b12}
B_{1,2}&\sim&\alpha^{|J_1|-|J_2|}\frac{I_1(\frac{t_1}{N},N)J_{C_2}(m_2)}{I_2(\frac{t_2}{N},N)J_{C_1}(m_1)}
\;.
\end{eqnarray}
Since the numbers $J_{C_i}(m_i),\;i=1,2$ are finite and positive, we have the following corollary of Theorem \ref{imalpha}.
\begin{cor}\label{practical1}
When the data belong to the open polytope  $C_i,\;i=1,2$, the Bayes factor $B_{1,2}$  is such that, 
when $\alpha\rightarrow 0$,  $$B_{1,2}\sim \alpha^{|J_1|-|J_2|}.$$
This implies in particular that, when the data is in both $C_i,\;i=1,2$  the Bayes factor always favours the sparser model.
\end{cor}
 
The proof follows immediately from (\ref{b12}). Moreover, when $\alpha\rightarrow 0$ and $|J_2|<|J_1|$, $B_{1,2}$ tends to $0$. This result has been well-known, at least numerically, for the class of decomposable models and in that case, it can be proved by expressing the Bayes factor as in (4.8) of \cite{mld09} and using the fact that $\Gamma(\alpha)\sim \alpha^{-1}$ as $\alpha\rightarrow 0$ (see Example 3 of Section 3 and Section \ref{sec:decomposable}). It has also been observed to hold numerically, most of the time, for hierarchical models. Computations illustrating the fact that the Bayes factor tends to favour the sparser models in the class of all hierarchical models can be found in \cite{mld09}, p. 3456. We have just shown that it actually always holds when the data is in  $C_1$ and in $C_2$. We will see in the next subsection that things are more delicate when the data belongs to the boundary of at least one of $\overline{C_1}$ or $\overline{C_2}$. 
\subsection{The case where the data belongs to a face of $\overline{C}_i,\; i=1,2$}

In this case, when $\alpha\rightarrow 0$, $\frac{\alpha m_i+t_i}{\alpha+N}$ converges to the boundary point $\frac{t_i}{N}$ of $C_i$ along the segment 
\begin{eqnarray}
\label{sa}
s(\alpha)=\frac{\alpha m_i+t_i}{\alpha+N}=\frac{\alpha}{\alpha +N}m_i+(1-\frac{\alpha}{\alpha +N})\frac{t_i}{N}.
\end{eqnarray}
 We need to study the limiting behaviour of $B_{1,2}$ when $\alpha\rightarrow 0$. To do so, we will use Theorem \ref{jmlambda2} to obtain the following  result.
\begin{theorem}\label{minimalpha}
Suppose that $\frac{t}{N}\in  \overline{C}\setminus C$ belongs to the relative interior of a face $F$ of dimension $k$. Then
\begin{equation}
\label{eq:minimalpha}
\lim _{\alpha\rightarrow 0}\alpha^{(|J|-k)} I(\frac{\alpha m+t}{\alpha+N},\alpha+N)
\end{equation}
 exists and is positive.
\end{theorem}
From Theorems \ref{imalpha} and \ref{minimalpha}, we immediately derive the following which is  the object of this paper. 
\begin{cor}\label{practical2}
Consider two hierarchical models $J_i,\;i=1,2$  of dimension $|J_i|$. Assume  that the data $\frac{t_i}{N}$  belongs to the relative interior of a face $F_i$ of $C_i$ of dimension $k_i$. Then the asymptotic behaviour of the Bayes factor $B_{1,2}$ when $\alpha\rightarrow 0$ is given by
\begin{eqnarray*}
B_{1,2}&\sim&D \alpha^{k_1-k_2}
\end{eqnarray*}
where $D$ is a finite positive constant. The Bayes factor favours the model which contains the data in the relative interior of the face of $C_i$ of smallest dimension.
\end{cor}
The proof is immediate.
According to Theorems \ref{imalpha} and \ref{minimalpha}, we have
\begin{eqnarray*}
B_{1,2}&=&\frac{I(m_2,\alpha)}{I(m_1,\alpha)}\frac{I(\frac{\alpha m_1+t_1}{\alpha+N},\alpha+N)}{I(\frac{\alpha m_2+t_2}{\alpha+N},\alpha+N)}
\sim \alpha^{|J_1|-|J_2|}\alpha^{(k_1-|J_1|)-(k_2-|J_2|)}=\alpha^{k_1-k_2},
\end{eqnarray*}
which proves our result.
\begin{remark}
We note that, if $\frac{t_i}{N}\in C_i,\;i=1,2,$ since $C_i$ is the face of $\overline{C_i}$ of dimension $J_i$, then $k_i=|J_i|$ and Corollary \ref{practical2} yields Corollary \ref{practical1}. For the same reason, Corollary \ref{practical2} also deals with the cases where  $\frac{t_i}{N}\in C_i$ for only $i=1$ or $i=2$.
\end{remark}

\begin{pff}\textsc{(of Theorem \ref{minimalpha})}
For simplicity of notation, we will write $n=|J|$ and  $y=t/N$ where  $N>0$ is fixed. Similarly, to simplify the expression of (\ref{sa}) above, we define
 $$\lambda=\frac{\alpha}{\alpha+N}\in (0,1)$$ 
 which implies  
 $\alpha=\frac{N\lambda}{1-\lambda}\;\;\mbox{ and}\;\; \alpha+N=\frac{N}{1-\lambda}.$ 
Our problem is then equivalent to showing that if $y$ belongs to the relative interior of $F$ and if $m\in C$ then 
$$\lim _{\lambda\rightarrow 0}\lambda ^{n-k} I(\lambda m+(1-\lambda)y,\frac{N}{1-\lambda})$$ 
exists and is positive.
 The idea of the proof is to consider the difference
 $$D(\lambda)=J_C(\lambda m+(1-\lambda)y)-(\frac{N}{1-\lambda})^{n}I(\lambda m+(1-\lambda)y,\frac{N}{1-\lambda})\;$$
 so that
  \begin{eqnarray}
  \label{diff} 
\hspace{5mm} (\frac{N}{1-\lambda})^{n}I(\lambda m+(1-\lambda)y,\frac{N}{1-\lambda})=J_C(\lambda m+(1-\lambda)y)-D(\lambda)\;.
\end{eqnarray}
Then, since from Theorem \ref{jmlambda2},  
$\lim_{\lambda\rightarrow 0} \lambda ^{n-k} J_C(\lambda m+(1-\lambda)y)$ exists and is positive, if we show that $\lim_{\lambda\rightarrow 0} \lambda ^{n-k} D(\lambda)$ exists and is  positive, Theorem \ref{minimalpha} will be proved. We proceed to do so now.

After the change of variable $\theta\in R^{n}\mapsto \frac{N}{1-\lambda}\theta\in R^{n}$ in $I(\lambda m+(1-\lambda)y,\frac{N}{1-\lambda})$, $D(\lambda)$ can be written
$$D(\lambda)= \int_{R^{|J|}}e^{\<\theta,\lambda m+(1-\lambda)y\>-h_C(\theta)}\left[1-\frac{e^{h_C(\theta)}}{L((1-\lambda)\theta/N)^{\frac{N}{1-\lambda}}}\right]d\theta$$
where $L(.)$ is defined in (\ref{lap}).

Consider the cone $A(f_k)=\{\theta\ ; \ h_C(\theta)=\<\theta,f_k\>\}.$ Let
$$D_k(\lambda)=\int_{A(f_k)}e^{\<\theta,\lambda m+(1-\lambda)y-f_k\>}\left[1-e^{\<\theta,f_k\>}(L((1-\lambda)\theta/N)^{-\frac{N}{1-\lambda}}\right]d\theta.$$ 
Since $D(\lambda)=\sum_{k\in {\cal E}}D_k(\lambda)$
we need only show that for each $k$, $\lim_{\lambda\rightarrow 0} \lambda ^{n-k}D_k(\lambda)$  exists.
Without loss of generality we can assume that  $f_k=0$ (if not, we replace $m$ and $y$ by $m-f_k$ and $y-f_k$). We want to prove that $\lim_{\lambda\rightarrow 0} \lambda ^{n-k}D_0(\lambda)$ exists and is finite. 
To do so, we split the cone $A(0)$ into a union of closed simplicial cones $(A(0)_s)_{s\in S_0}$ with disjoint interiors so that $A(0)=\cup_{s\in S_0}A(0)_s$ and we write $D_0(\lambda)=\sum_{s\in S_0}\Delta_s(\lambda)$ where 
\be\label{AA's}
\Delta_s(\lambda)=\int_{A(0)_s}e^{\<\theta,\lambda m+(1-\lambda)y\>}\left[1-(L((1-\lambda)\theta/N)^{-\frac{N}{1-\lambda}}\right]d\theta.
\ee 
We need to prove that
 $\lim_{\lambda\rightarrow 0} \lambda ^{n-k}\Delta_s(\lambda)$ exists and is non negative. 
The simplicial cone $A(0)_s$ is defined by $n$ independent linear vectors $g_j$ in $R^{n}$ as 
$$A(0)_s=\{\theta=\lambda_1g_1+\cdots+\lambda_n g_n \in R^n\ ;\ \lambda_1,\ldots,\lambda_n\leq 0\}.$$ Without loss of generality we assume $|\det[g_1,\ldots,g_n]|=1.$  For simplicity we write 
$$m_j(\lambda)=\<g_j,\lambda m+(1-\lambda)y\>.$$ 
We observe that $m_j(\lambda)\geq 0$ on $A(0)_s^o\supset A(0)^o.$ 
Thus by the change of variable $\theta=\lambda_1g_1+\cdots+\lambda_n g_n\mapsto (\lambda_1,\ldots,\lambda_n)$ (giving $d\theta=d\lambda_1\ldots d\lambda_n$ ) we obtain 
$$\int_{A(0)_s}e^{\<\theta,\lambda m+(1-\lambda)y\>}d\theta=\frac{1}{\prod_{j=1}^n m_j(\lambda)}.$$  
Now, in the integral $\int_{A(0)_s}e^{\<\theta,\lambda m+(1-\lambda)y\>}(L((1-\lambda)\theta/N)^{-\frac{N}{1-\lambda}}d\theta$ in the right-hand side of (\ref{AA's}), we make the further change of variable  $$(\lambda_j,\;j=1,\ldots,n)\mapsto (u_j=e^{(1-\lambda)\lambda_j/N},\;j=1,\ldots,n)$$ 
giving $d\lambda_1\ldots d\lambda_n=(\frac{N}{1-\lambda})^n \frac{du_1\ldots du_n}{u_1\ldots u_n}$ with $u_j \in (0,1).$ 
Moreover, since $L(\theta)=1+\sum_{i\not =0,i\in {\cal I}}e^{\<\theta,f_i\>}$, we have   $$L(\frac{1-\lambda}{N}\theta)=L((\log u_1) g_1+\cdots + (\log u_n) g_n)=1+\sum_{i\neq 0}\prod_{j=1}^nu_j^{a_{ij}}$$ where $a_{ij}=\<g_j,f_i\>=\<g_j,f_j+f_i-f_j\>\geq 0$ since $\langle g_j,f_j-x\rangle=0$ is a supporting hyperplane of $C$ at $f_j$. To simplify notation, write
$$h_{\lambda}(u_1,\ldots,u_n)=(1+\sum_{i\neq 0}\prod_{j=1}^nu_j^{a_{ij}})^{-\frac{N}{1-\lambda}}.$$
Recalling that $\alpha+N=\frac{N}{1-\lambda}$ we then have 
\begin{equation}
\label{alexa}
\Delta_s(\lambda)= \frac{1}{\prod_{j=1}^n m_j(\lambda)}\left[1-K(\lambda)\right]
\end{equation}
 where 
 \begin{eqnarray*}
 K(\lambda) &=& (\alpha+N)^n \prod_{j=1}^n m_j(\lambda)\times \\
 \int_0^{1}\ldots \int_0^{1}&&u_1^{(\alpha+N) m_1(\lambda)-1}\ldots u_n^{(\alpha+N) m_n(\lambda)-1} h_{\lambda}(u_1,\ldots,u_n)du_1\ldots du_n\;.  
 \end{eqnarray*}
 Note that $h_{\lambda}<1$ and that $K(\lambda)$ can be seen as the expectation of $h_{\lambda}(U_1,\ldots,U_n)$ where $U_1,\ldots,U_n$ are independent random variables with density  $f(u_j)=(\alpha+N)m_j(\lambda)u_j^{(\alpha+N)m_j(\lambda)-1}$ on $(0,1).$
 
Now to determine the behaviour of $\Delta_s(\lambda)$ when $\lambda\rightarrow 0$, we recall that since $A(0)_s$ is simplicial, the supporting hyperplane of $A(0)_s^0$ are the hyperplanes defined by 
  \be
  \label{plane}\{x\in {R}^n; \<g_j,x\>=0\},\;\;j=1,\ldots,n.\ee
  The data point $y$ either belongs to a face of $A(0)_s^0$ of dimension less than $n$ or it belongs to its interior. If it belongs to its interior, then for all 
  $j=1,\ldots,n$, 
  $$m_j(\lambda)\rightarrow_{\lambda\rightarrow 0}\<g_j,y\>$$
  and standard reasoning shows that the limit of $K(\lambda)$ exists and is equal to $K=E(h_{\lambda}(u_1,\ldots,u_n))$ where $U_j,\;j=1,\ldots,n$ follow independent distributions with density $f_j(u_j)=N\<g_j,y\>u_j^{N\<g_j,y\>-1}$ and therefore 
  $$\mbox{lim}_{\lambda\rightarrow 0}\lambda^{n-k}\Delta_s(\lambda)\rightarrow 0.$$
If $y$ belongs to a face $F^0$ of $A(0)_s^0$, the dimension $k^0$ of $F^0$ is greater than or equal to $k$ so that $n-k^0\leq n-k$. Therefore $F^0$ is contained in the intersection of $n-k^0$ hyperplanes of the type (\ref{plane}).  
 Without loss of generality we can assume that these supporting hyperplanes of $A(0)^o_s$ have been numbered so that the first $n-k^0$ are those  containing $y$, that is,
  $$\{x\in {R}^n; \<g_j,x\>=0\}$$
 for $j=1,\ldots,n-k^0.$ As a consequence we have that $\lim_{\lambda\rightarrow 0}m_j(\lambda)=0$ 
 for $j=1,\ldots,n-k^0$ and $\lim_{\lambda\rightarrow 0}m_j(\lambda)=\<g_j, y\> >0$ if $j=n-k^0+1,\ldots,n.$ Thus the limiting distribution of $U_j$ when $\lambda\rightarrow 0$ is the Dirac mass at 0 if $j\leq n-k^0$ and is the distribution with density $N\<g_j, y\>u_j^{N\<g_j, y\>-1}$ on $(0,1)$ if $n-k^0<j\leq n.$ It is easy to show that 
 \begin{equation}\label{DK}K=\lim_{\lambda\rightarrow 0}K(\lambda)\end{equation}
  exists and is
 $$N^{k^0}\prod_{j=n-k^0+1}^n\<g_j, y\>\int_0^1\ldots\int_0^1h_0(0,\ldots,0,u_{n-k^0+1},\ldots,u_n)\prod_{j=n-k^0+1}^nu_j^{N\<g_j, y\>-1}du_j.$$ 
 Recall that for $1\leq j\leq n-k^0$ we have $m_j(\lambda)=\lambda\<g_j,m\>.$ We get from (\ref{alexa}) that $$\lim_{\lambda\rightarrow 0}\lambda^{n-k^0}\Delta_s(\lambda)=\frac{1}{\prod_{j=1}^{n-k^0}\<g_j,m\>\prod_{j=n-k^0+1}^n\<g_j,y\>}(1-K).$$
 Since $h_0\leq 1$, this limit exists and is nonnegative and so does \newline $\lim_{\lambda\rightarrow 0}\lambda^{n-k}\Delta_s(\lambda)$. Moreover, since $I(m,\alpha)$ is always positive,  we also have that
 $\lim_{\alpha\rightarrow 0}\alpha^{n-k} I(\frac{\alpha m+t}{\alpha+N},\alpha+N)$ always exists and is nonnegative.
 
 We now need to prove that the latter is actually positive. To do so, it is sufficient to prove that $K$ in (\ref{DK}) is strictly less than $1$ for the $K(\lambda)$ corresponding to  at least one of the simplicial cones $A(f)_s$.
  To do so, let us first remark that if $y$ coincides with a vertex $f$ of $C$, then $\mbox{lim}_{\lambda\mapsto 0}m_j(\lambda)=0$ for all $j=1,\ldots,n$ so that, as $\lambda$ tends to $0$, all densities $(\alpha+N)m_j(\lambda)u_j^{(\alpha+N)m_j(\lambda)-1}{\bf{1}}_{(0,1)}d u_j$ tend to the Dirac mass at $0$ and $K(\lambda)$ tends to 1. So, for such an $A(f)_s$, $K=1$. Clearly there exists an $f_{i_0}$ such that $y$ does not coincide with $f_{i_0}$ and also such that $k^0=k$ (otherwise, $F$ would not contain $y$ in its relative interior). In such a case, the number $n-k$ of faces of $A(f_{i_0})_s^0$ containing $y$ is strictly less than $n$ and for $j=1,\ldots,n-k$,
  $\langle g_j,f_{i_0}\rangle=0$ and therefore
  $$1+\sum_{i\not =0}\prod_{j=1}^nu_j^{a_{ij}}\geq 1+\prod_{j=n-k+1}^nu_j^{a_{i_oj}}>1$$ 
  and $h_{\lambda}(u_1,\ldots,u_n)<1.$
  Since  the densities of $u_j,\;j=n-k+1,\ldots,n$ are  proper Beta densities when $\lambda\rightarrow 0$, the limit $K$ is strictly less than 1 and we have now proved that $\lim _{\alpha\rightarrow 0}\alpha^{n-k} I(\frac{\alpha m+t}{\alpha+N},\alpha+N)$ is strictly positive. 
 \end{pff}

\subsection{The results of Steck and Jaakola \cite{sj02} as a particular case }
\label{steck}
In \cite{sj02} Steck and Jaakola study the behaviour of the Bayes factor for two Bayesian network models differing by one edge only, when $\alpha\rightarrow 0$. They show it is equivalent to the problem of comparing two Bayesian network models with three variables indexed by $\{a,b,c\}$. The first model has directed edges $(b,a),\;(b,c)$ and $(a,c)$. The second model has directed edges $(b,a)$ and $(b,c)$. These two Bayesian network models are Markov equivalent to the two hierarchical (in fact graphical) models $J_1$ and $J_2$ with, respectively, generating sets ${\cal D}_1=\{abc\}$ and ${\cal D}_2=\{ab, bc\}$. Moreover on these two models, the prior in \cite{sj02} is equivalent to ours. We must then be able to compare their result given in  Proposition 1 of \cite{sj02} and our result given in Corollary \ref{practical2}. 
To give their results Steck and Jaakola \cite{sj02} introduce the quantity
\begin{equation}
d_{EDF}=\sum_{i\in {\cal I}}\delta(n(i))-\sum_{i_{ab}\in {\cal I}_{ab}}\delta(n(i_{ab}))-\sum_{i_{bc}\in {\cal I}_{bc}}\delta(n_{i_{bc}})+\sum_{i_b\in {\cal I}_b}\delta(n(i_b))
\end{equation}
where $\delta(.)$ is an indicator function which is such that $\delta(x)=0$ if $x=0$ and $\delta(x)=1$ otherwise.
They state that the Bayes factor $B_{1,2}$ behaves as follows
$$\mbox{lim}_{\alpha\rightarrow 0}B_{1,2}=\left\{\begin{array}{cc}0&\mbox{if}\;\;d_{EDF}>0\\+\infty&\mbox{if}\;\;d_{EDF}<0\end{array}
\right .$$
This result coincides with our Corollaries \ref{practical1} and \ref{practical2} for three variable models. In fact, we are going to show the following.
\begin{prop}
\label{EDF}
Consider the two decomposable graphical models on three variables, $J_1$ and $J_2$, as defined above. If the data belongs to faces of dimension $k_1$ and $k_2$ of, respectively,  $C_1$ and $C_2$, then we have
$$d_{EDF}=k_1-k_2.$$
\end{prop}

\begin{pff}
 The Bayes factor is equal to
\begin{eqnarray*}
\frac{I(\frac{\alpha m_1+t_1}{\alpha+N},\alpha+N)I(m_2,\alpha)}{I(m_1,\alpha)I(\frac{\alpha m_2+t_2}{\alpha+N},\alpha+N)}
\end{eqnarray*}
where the form of the normalizing constants $I(m,\alpha)$ for decomposable models is well-known (see for example equation (4.8) of \cite{mld09}).
When $\alpha\rightarrow 0$, from Theorem \ref{imalpha}, we know that
\begin{eqnarray*}
\frac{I(m_2,\alpha)}{I(m_1,\alpha)}\sim \alpha^{|J_1|-|J_2|}
\end{eqnarray*}

Expressed in terms of cell counts for the full table, for the $b$-, $ab$- and $bc$- marginal tables, we have
\begin{equation}\label{lots}
\frac{I(\frac{\alpha m_1+t_1}{\alpha+N},\alpha+N)}{I(\frac{\alpha m_2+t_2}{\alpha+N},\alpha+N)}
=\frac{\prod_{i\in {\cal I}}\Gamma(\alpha m(i)+n(i))\prod_{i_b\in {\cal I}_b}\Gamma(\alpha m(i_b)+n(i_b))}{\prod_{i_{ab}\in {\cal I}_{ab}}\Gamma(\alpha m(i_{ab})+n(i_{ab}))\prod_{i_{bc}\in {\cal I}_{bc}}\Gamma(\alpha m(i_{bc})+n(i_{bc}))}
\end{equation}
If for some $D=\emptyset, ab, bc, b$, the marginal cell count $n(i_{D})$ is different from 0, when $\alpha\rightarrow 0$, 
$\Gamma(\alpha m(i_D)+n(i_D))\rightarrow \Gamma(n(i_D))$ which is finite. If $n(i_D)= 0$, then $\Gamma(\alpha m(i_D)+n(i_D))\sim \frac{1}{\alpha m(i_D)}$. It follows from (\ref{lots}) that, when $\alpha\rightarrow 0$,
$B_{1,2}\sim \alpha^{q}$ 
where 
\begin{eqnarray*}
q&=&[|J_1|-\sum_{i\in {\cal I}}(1-\delta(n(i)))]\\
&&\hspace{.5cm} -[|J_2| -\sum_{i\in {\cal I}_{ab}}(1-\delta(n(i_{ab})))-\sum_{i_{bc}\in {\cal I}_{bc}}(1-\delta(n(i_{bc})))+\sum_{i_{b}\in {\cal I}_{b}}(1-\delta(n(i_{b}))].
\end{eqnarray*}
Let $C_i,\; i=1,2$ be the interior of the  convex hull corresponding to model $J_i$. 
Consider model $J_1$ first. It is immediate to see that, following the notation of (\ref{F0}) and (\ref{Fj}) in Section 5 below
\begin{eqnarray*}
n(000)&=&g_{0,C_1}\\
n(i)&=&g_{i,C_1}, \;i\in{\cal I}
\end{eqnarray*} 
and according to Theorem \ref{generalfacets} $n(000)=0$ and $n(i)=0$ are the equations of the facets of the polytope $C_1$. Therefore the dimension of the space minus the number of distinct facets the data belongs to, is equal to the dimension of the face of $\overline{C_1}$ containing the data, that is,
\begin{equation}
\label{k1}
|J_1|-\sum_{i\in {\cal I}}(1-\delta(n(i)))=\sum_{i\in {\cal I}}\delta(n(i))=k_1.
\end{equation}
Similarly, for model $J_2$, according to Theorem \ref{generalfacets}, the equations of the facets of $\overline{C_2}$ are given by 
$$n(i_{ab})=0,\;i_{ab}\in {\cal I}_{ab}\;\;\mbox{and}\;\;n(i_{bc})=0,\;i_{bc}\in {\cal I}_{bc}.$$
The facets containing the data are therefore those  defined by $n(i_{ab})=0$ or $n(i_{bc})=0$. This does not mean, however, that 
$$|J_2|-(1-\sum_{i_{ab}\in {\cal I}_{ab}}\delta(n(i_{ab})))-\sum_{i_{bc}\in {\cal I}_{bc}}(1-\delta(n(i_{bc})))$$
represents the dimension of the face containing the data. Indeed, if for some $i^0_b\in {\cal I}_b$, we have $n(i^0_b)=0$, this means that $n(i_{ab})=0$ also whenever $i_{b}=i^0_b$ and also $n(i_{bc})=0$ whenever $i_b=i^0_b$. Then clearly, one of the equations $n(i_{ab})=0$ or $n(i_{bc})=0$ is redundant and we subtract $1-\delta(n(i^0_b))$ for the count of facets defining the position of the data. It is clear then that
$$|J_2| -\sum_{i\in {\cal I}_{ab}}(1-\delta(n(i_{ab}))-\sum_{i_{bc}\in {\cal I}_{bc}}(1-\delta(n(i_{bc}))+\sum_{i_{b}\in {\cal I}_{b}}(1-\delta(n(i_{b}))=k_2,$$
which, together with (\ref{k1}) proves the proposition.
\end{pff}

In fact Proposition \ref{EDF} can be extended to the following  general result. Let ${\cal C}_i$ and ${\cal S}_i$ the set of cliques and separators of the decomposable model $J_i,\;i=1,2$. We define the effective degrees of freedom to be the following sum $d_{EDF}:$ 
\begin{eqnarray*}
d_{EDF}&=&\sum_{C\in {\cal C}_1}\sum_{i_{C}\in {\cal I}_{C}}\delta(n(i_{C}))-\sum_{S\in {\cal S}_1}\sum_{i_{S}\in {\cal I}_{S}}\delta(n(i_{S}))\\
&&\hspace{1cm}-\Big(\sum_{C\in {\cal C}_2}\sum_{i_{C}\in {\cal I}_{C}}\delta(n(i_{C}))-\sum_{S\in {\cal S}_2}\sum_{i_{S}\in {\cal I}_{S}}\delta(n(i_{S}))\Big).
\end{eqnarray*}
\begin{prop}
\label{generalEDF}
Consider two arbitrary decomposable graphical models $J_1$ and $J_2$ such that the data belongs to faces of dimension $k_1$ and $k_2$  of $C_1$ and $C_2$ respectively. 
Then, the following relation holds:
$$d_{EDF}=k_1-k_2.$$
\end{prop}
The proof of this proposition follows parallel lines to the proof given above for the two particular models given in \cite{sj02}. We therefore have a quick and easy way to know the behaviour of the Bayes factor between two decomposable models.

\section{Facets of $\overline{C}$ for some hierarchical models}
We now turn our attention to the identification of the facets of $\overline{C}$. Knowing the facets of $\overline{C}$ is crucial since faces are intersection of facets. Facets of $\overline{C}$ have been much studied by geometers and in Section 5.3, we will recall some known results on these facets when the model is binary and governed by a cycle of order $n\geq 3$. But before doing so, we give two new results on facets of polytopes associated to our models. In Theorem \ref{generalfacets}, we 
 identify a category of facets which is common to all discrete hierarchical models. In  Corollary \ref{facedecomp}, we show that for decomposable graphical models, the only facets of $C$ are given by the category of facets given in Theorem \ref{generalfacets}.
 \subsection{Facets common to all hierarchical models}
 Let $\mathcal{D}$ the set of subsets of $V$ defining the hierarchical model. Let $\mathcal{A}$ be the family of maximal elements of $\mathcal{D}$. For the subclass of graphical models Markov with respect to a graph $G$, $\mathcal{A}$ is the set of cliques of $G$. This set is traditionally denoted ${\cal C}$ but in this particular subsection, to avoid confusion between a clique $C\in {\cal C}$ and the polytope $C$, we use the notation $A\in {\cal A}$.
 
 \noindent For each $D\in \mathcal{D}$ and each $j_0\in J$ such that $S(j)\subset D$ define the  affine forms 
\begin{eqnarray} 
\label{F0}g_{0,D}(m)&=&1+\sum_{ j; S(j)\subset D}(-1)^{|S(j)|}m_j\\
\label{Fj}
g_{j_0,D}(m)&=&\sum_{j; S(j)\subset D,\ j_0\tl j}(-1)^{|S(j)|-|S(j_0)|}m_j
\end{eqnarray}
 Of course,  for $j_0\in J$ the form $g_{j_0,D}$ is not only affine, but is also linear. 
 In this subsection, we will use $g_{j,A}$ only for $A\in \mathcal{A}$ but as we shall see in Subsection \ref{sec:decomposable}, $g_{j,S}$ when $S$ is a minimal separator play an important role also even though $S\not \in {\cal A}$. In the next theorem, we consider the following affine hyperplanes of $R^J$ 
 $$H(j,A)=\{m\in R^J\; ; \; g_{j,A}(m)=0\},\;j\in J\cup\{0\},\;A\in \mathcal{A}$$ and we prove that 
 \begin{equation}
 \label{fj}
 F(j,A)=H(j,A)\cap \overline{C}
 \end{equation}
  is a facet of the convex set  $\overline{C}$ with extreme points $f_i=\sum_{j\tl i}e_j.$ Recall that for $T\subset V$ the index set $I_T$ means $\prod_{v\in T} I_v.$
 
\begin{theorem} \label{generalfacets}
Let $A\in \mathcal{A}$ be the set of maximal elements of ${\cal D}$ defining a general hierarchical model. Let $j_0\in J\cup\{0\}$ such that $ S(j_0)\subset A$ and let  $i\in I.$ Then $g_{j_0,A}(f_i)$ can only take values 0 or 1. More precisely, the following holds:
 \begin{enumerate}
 \item[1.] $g_{j_0,A}(f_i)=1$ if and only if $j_0\tl i$ and $S(i)\cap A= S(j_0);$  
 \item[2.] there are exactly $|I|-|I_{V\setminus C}|$ vectors  $f_i$'s such that $g_{j_0,A}(f_i)=0$; 
 \item[3.]  the set $F(j_0,A)$ as defined in (\ref{fj}) is a facet of the polytope $\overline{C}.$ 
 \end{enumerate}
 \end{theorem}
 

 
 
Before giving the proof of this theorem, let us illustrate its results. We consider a simple decomposable model and list the various faces and the $f_i$'s that belong to them.
 
\vspace{4mm}\noindent {\it Example:} the $A_3$ graph. Consider  $\stackrel{a}{\bullet}-\stackrel{b}{\bullet}-\stackrel{c}{\bullet}$ and assume that we are in the binary case $I_a=I_b=I_c=\{0,1\}.$ For each $F(j,A)$, we are going to list the $f_i$'s that belong to it.
In this example, $I$ is identified with the power set of $V=\{a,b,c\}$ and $J$ is identified with the set of nonempty complete subsets $D$ of $A_3$, namely $a,b,c,ab,bc.$ 
In a five-dimensional space with basis $e_a,e_b,e_c,e_{ab},e_{bc}$ the eight vectors $f_T$ are 
$$0,\ f_a=e_a,\ f_b=e_b,\ f_c=e_c,\ $$$$f_{ab}=e_a+e_b+e_{ab},\ f_{ac}=e_a+e_c,\ f_{bc}=e_b+e_c+e_{bc},\ f_{abc}=e_a+e_b+e_c+e_{ab}+e_{bc}.$$ 
Since we have two cliques of size 2 the number of facets $F_{D,C}$ is $2^2+2^2=8.$ They are described as follows (we adopt the following short notation : $F^*_{D,A}$ is the set of $T$'s contained in $V$ such that $f_T\in F_{D,A}).$

\begin{eqnarray*}
F^*_{\emptyset,ab}&=&\{a,b,ab,bc,ac,abc\},\ 
F^*_{a,ab}=\{\emptyset, b,c,ab,bc,abc\},\\  
F^*_{b,ab}&=&\{\emptyset,a,c,ab,ac,abc\},\  
F^*_{ab,ab}=\{\emptyset,a,b,c,bc,ac\},\\ 
F^*_{\emptyset,bc}&=&\{b,c,ab,bc,ac,abc\},\  
F^*_{c,bc}=\{\emptyset,a,b,ab,bc,abc\},\\  
F^*_{b,bc}&=&\{\emptyset,a,c,bc,ac,abc\},\ 
F^*_{bc,bc}=\{\emptyset,a,b,c,ab,ac\}.
\end{eqnarray*}

\begin{pff} \textsc{of Theorem \ref{generalfacets} } The proof is long and we will only give it for $j_0\in J.$ We skip the case $j_0=0$ since it is entirely analogous. We have
\begin{eqnarray*}
g_{j_0,A}(f_i)&=&\sum_{j_0\tl j\atop S(j)\subset A}(-1)^{|S(j)|-|S(j_0)|}\langle e_{j},\sum_{j^{'}\tl i}e_{j^{'}}\rangle\\
&=& \sum_{j_0\tl j\tl i\atop S(j)\subset  A}(-1)^{|S(j)|-|S(j_0)|}\;.
\end{eqnarray*}
From (\ref{TL}) the set of $j\in J$ such that $j_0\tl j\tl i$ is non empty if and only if $j_0\tl  i.$ Thus $g_{j_0,A}(f_i)=0$ if $j_0\tl  i$ is false. Suppose now that $j_0\tl  i.$ Thus $S(j_0)\subset A\cap S(i).$ 
If $S(i)\cap A=S(j_0)$, then the only $j$ satisfying the conditions of the sum above is $j=j_0$ and clearly $g_{j_0,A}(f_i)=1.$ If $S(i)\cap A\not =S(j_0)$, then  any $j$ such that $j_0\tl j\tl  i$ and $S(j)\subset S(i)\cap A$  can be written $j=(i_{S(j_0)},i_{S(j)\setminus S(j_0)},0)$. This implies that the number $\sum_{j:j_0\tl j\tl i\atop S(j)\subset  A}(-1)^{|S(j)|-|S(j_0)|}$ can be computed by the principle of inclusion exclusion and is equal to zero. We have just proved part 1. of the theorem, that is 
\begin{equation}
\label{a}
g_{j_0,A}(f_i)=\left\{\begin{array}{cc}1&\mbox{if}\;\;j_0\triangleleft i\;\;\mbox{and}\;\;S(i)\cap A=S(j_0)\\
0&\mbox{otherwise}\;. 
\end{array}\right.
\end{equation} 
\vspace{3mm}

\noindent From Part 1., the number of $f_i$ contained in $H(j_0,A)$ is equal to $|I|$ minus the number of $f_i$ such that  $g_{j_0,A}(f_i) =1$. Since $S(i)\cap A=S(j_0)$ and $j_0\triangleleft i$, such $i$'s are identified by $i_{V\setminus A}$. Clearly their number is equal to $|I_{V\setminus A}|$. This proves Part 2.
\vspace{3mm}

\noindent Let us now prove Part 3.
We know from Proposition \ref{extreme} that the $f_i$ are the extreme points of $\overline{C}$. 
Therefore, since  $g_{j_0,A}(f_i)\geq 0$ for all $f_i$ and  $g_{j_0,A}(f_i)=0$ for some $f_i$, then $H(j_0,A)$ is a supporting hyperplane of $\overline{C}.$

The more delicate part of the theorem is to show that $F(j_0,A)$ is a facet of $\overline{C}$. This is equivalent to saying  that if $j_{0}\in J$ and $A\in \mathcal{A}$ are such that $S(j_0)\subset A$ then $H(j_0,A)$ contains enough points $f_i$ which affinely generate it. Since $j_0\in J$  is not zero, then $f_0=0$ is in $H(j_0,A)$ which is therefore a linear space. To prove that $F(j_0,A)$ is a facet of $\overline{C}$, we want to prove that it linearly generates $H(j_0,A)$. This is equivalent to proving the following statement. 
 
 \textit{Statement S:} \textit{If $h\in H(j_0,A)$  is orthogonal to all elements $f_i\in H(j_0,A)$ then $h=0.$ }

 We write $h=\sum_{j\in J}h_je_j.$ We prove  that $h_j=0$ for all $j\in J$ in three steps.

 \vspace{4mm}\noindent \textsc{Step 1.} We prove that if $h_j\neq 0$ then $j_0\tl j.$ Let $\mathcal{N}=\{j\in J\ ;\ h_j\neq 0\}.$ Let $j_1$ be a minimal element of $\mathcal{N}$ (we mean that $h_{j_1}\neq 0$ and that $h_j=0$ for all $j\tl j_1$ with $j\neq j_1).$ Therefore 
 $$\<f_{j_1},h\>=\sum_{j\tl j_1}h_j
= h_{j_1}\neq 0.$$ Since $h$ is orthogonal to all $f_i\in H(j_0,A)$ we get that $f_{j_1}$ is not in $H(j_0,A)$ and, as we showed earlier in this proof, this implies that $j_0\tl j_1$. Now let $j$ such that $h_j\neq 0.$ There exists necessarily a minimal element $j_1$ of $\mathcal{N}$ such that $j_1\tl j.$ Therefore $j_0\tl j_1 \tl j$ and Step 1 is proved.

 \vspace{4mm}\noindent \textsc{Step 2.} We prove that if $j_0\tl j$ and $S(j)\subset A$ we have $h_j=0.$  Let
 $\varphi(j)=\sum_{j_0\tl j'\tl j}h_{j'}.$ If $j\neq j_0$ we have the following equalities
 $$\varphi(j)=\sum_{j_0\tl j'\tl j}h_{j'}\stackrel{(1)}{=}\sum_{j'\tl j}h_{j'}\stackrel{(2)}{=}\<f_{j},h\>\stackrel{(3)}{=}0.$$ 
 Indeed, (1) is a consequence of Step 1, (2) is by definition of $f_j$. For (3), we see that since
 $j\neq j_0$,  $S(j)\cap A\not =S(j_0)$ and therefore by (\ref{a})), $f_j\in H(j_0,A)$. Since $h$ is orthogonal to any element of $H(j_0,A)$, (3) follows. However if $j= j_0$ then $\varphi(j_0)=h_{j_0}.$ The inclusion exclusion principle applied to $\varphi(j)$ yields, for $j_0\tl j$ and $S(j)\subset A$, 
 \be\label{STWO}h_j=\sum_{j_0\tl j'\tl j}(-1)^{|S(j)|-|S(j')|}\varphi(j')=(-1)^{|S(j)|-|S(j_0)|}h_{j_0}.\ee We now use the hypothesis
 $h\in H(j_0,A)$ that is 
 $$0=\<g_{j_0,A},h\>=\sum_{j_0\tl j, \ S(j)\subset A}(-1)^{|S(j)|-|S(j_0)|}h_j=h_{j_0}\sum_{j_0\tl j, \ S(j)\subset A}1$$
 As a consequence $h_{j_0}=0$ and (\ref{STWO}) gives Step 2. 
 
 \vspace{4mm}\noindent \textsc{Step 3.} We prove that if $j_0\tl j$ and $S(j)\not \subset A$ we have $h_j=0.$ Once we prove this, Statement S will be shown and this will complete the proof of Theorem \ref{generalfacets}. We prove Step 3 by induction on the size $k$ of the set $S(j)\setminus A.$ 
 
 \noindent For $k=0$, it is Step 2. To understand the principle of the proof it is wise to give this proof first for $k=1.$  Although throughout the paper the symbol $i'\tl i$ was used only when $i'\in J$ it makes sense even if $i'$ and $i$ are in $I$ and we can write
 $$\<f_i,h\>=\sum_{i'\tl i}h_{i'}$$ with the  convention 
 \begin{equation}
 \label{c}
 h_{i'}=0\;\;\mbox{for}\;\; i'\in I\setminus J.
 \end{equation} 
 We fix now $i\in I$ such that $j_0\tl i$ and that $S(i)=A\cup\{v\}$ where  $v\in V\setminus A.$ For $S_1\subset A\setminus S(j_0)$ consider the unique $i(S_1)$ such that $j_0\tl i(S_1)\tl i$ and  $S(i(S_1))=S(j_0)\cup S_1\cup\{v\}.$
 Define now 
 $$\varphi_v(S_1)=\<f_{i(S_1)},h\>\stackrel{(1)}{=}\sum_{j_0\tl i'\tl i(S_1)}h_{i'}
 \stackrel{(2)}{=}\sum_{S\subset S_1}h_{i_{S(j_0)\cup S\cup\{v\}}}.$$
 In this equality, (1) follows from Step 1 and (\ref{c}) while (2) follows from two remarks. The first remark is that if $j_0\tl i'\tl i(S_1)$ then $i'$ is entirely determined by its support $S(i')$ since $i'\tl i.$ This support has two possible forms: either $S(i')=S(j_0)\cup S$  or $S(i')=S(j_0)\cup S\cup\{v\},$ with $S\subset S_1.$ The second remark is that if $i'$ has a support of the form $S(i')=S(j_0)\cup S$ then $h_{i'}=0$. This follows from Step 2 if $i\in J$ and  from (\ref{c}) if $i\not\in J$.
 From (\ref{a}), since $j_0\tl i(S_1)$ we have $\varphi_v(S_1)=\<f_{i(S_1)},h\>=0$ if and only if $S(i(S_1))\cap A\neq S(j_0),$ that is if and only if $S_1\neq \emptyset.$ Moreover,  $\varphi_v(\emptyset)=h_{j_{S(j_0)\cup\{v\}}}.$ The inclusion exclusion principle applied to $\varphi_v(S)$ therefore implies that 
 $$h_{j_{S(j_0)\cup S\cup\{v\}}}=(-1)^{|S|}h_{j_{S(j_0)\cup\{v\}}}$$ for all $S\subset A\setminus S(j_0).$ 
 We now  apply this last equality to $S=A\setminus S(j_0)$ itself. Because of the maximality of $A\in \mathcal{A}$ the set  $A\cup\{v\}$ is not in $\mathcal{D}$ and therefore $h_{j_{C\cup\{v\}}}=0$ from our convention. As a consequence $h_{j_{S(j_0)\cup\{v\}}}=0$  and also $h_j=0$ for all $j_0\tl j\tl i.$

 This settles the case where $S(i)=A\cup\{v\}$ where  $v\in V\setminus A.$ 
 We now make the following induction hypothesis on $k$: if $i\in I$ is such that $A\subset S(i)$ and such that $S(i)\setminus A$ has $k$ elements then $h_j=0$ for all $j\in J$ such that $j_0\tl j\tl i.$  We assume that this induction hypothesis is true up to $k-1.$ We denote 
 $$ S(i)=A\cup\{v_1,\ldots,v_k\}.$$  
 For $S_1\subset A\setminus S(j_0)$ consider $i(S_1)$ such that $j_0\tl i(S_1)\tl i$ and defined by $S(i(S_1))=S(j_0)\cup S_1\cup \{v_1,\ldots,v_k\}.$ We now fix a subset $R$ of $\{v_1,\ldots,v_k\}.$
 Define $$\varphi_R(S_1)=\sum_{S\subset S_1}h_{i_{S(j_0)\cup S\cup R}}.$$
 We have
 $$\<f_{i(S_1)},h\>=\sum_{j_0\tl i'\tl i(S_1)}h_{i'}=\sum_{R\subset \{v_1,\ldots,v_k\}} \varphi_R(S_1)$$
 Now the induction hypothesis implies that $\varphi_R(S_1)=0$ if $|R|<k.$ We get that 
 $\<f_{i(S_1)},h\>=\varphi_{\{v_1,\ldots,v_k\}}(S_1).$ Since $j_0\tl i(S_1)$ we have $\varphi_{\{v_1,\ldots,v_k\}}(S_1)=\<f_{i(S_1)},h\>=0$ if and only $S(i(S_1))\cap A\neq S(j_0),$ that is if and only if $S_1\neq \emptyset.$  Similarly to the case $k=1$ we have $\varphi_{\{v_1,\ldots,v_k\}}(\emptyset)=h_{j_{S(j_0)\cup\{v_1,\ldots,v_k\}}}.$ The inclusion exclusion principle therefore implies that 
 $$h_{j_{S(j_0)\cup S\cup\{v_1,\ldots,v_k\}}}=(-1)^{|S|}h_{j_{S(j_0)\cup\{v_1,\ldots,v_k\}}}$$ 
 for all $S\subset A\setminus S(j_0).$ Again $h_{j_{A\cup\{v_1,\ldots,v_k\}}}=0$ since $C\cup\{v,1,\ldots,v_k\}$ is not in $\mathcal{D}$ 
 and this leads to $h_j=0$ if $j_0\tl j\tl i.$ The induction is extended and Statement S is proved as well as Theorem \ref{generalfacets}.  
 \end{pff}
 
 \subsection{Facets of $\overline{C}$ when $G$ is decomposable}
\label{sec:decomposable}
When the graph $G$ is decomposable, the normalizing constant $I(m,\alpha)$ is the normalizing constant of the hyper Dirichlet as defined in \cite{dl93}.
In the theorem below, we restate, in our present notation, the expression of $I(m,\alpha)$ as given in \cite{mld09}, Formula (4.8) and directly derive the form of $J_C(m)$ for decomposable models. A corollary giving the facets of $\overline{C}$ when the model is decomposable follows immediately from the theorem.

 \begin{theorem} \label{MLD}Let $(V,\mathcal{E})$ be a decomposable graph, let $\mathcal{D}$ be the family of the complete subsets of $V,$  let $\mathcal{C}$ be the family of its cliques, let $\mathcal{S}$ be the family of its minimal separators and let $\nu(S)$ be the multiplicity of the minimal separator $S.$ 
 Then for $m$ in the interior $C$ of the convex hull of the $f_i$'s we have 
 \begin{eqnarray}
 \hspace{2cm}I(m,\alpha)&=&\int_{R^J}e^{\alpha\<\theta,m)}L(\theta)^{-\alpha}d\theta\nonumber\\
 &=&\frac{\prod_{C\in \mathcal{C}} \Gamma(\alpha g_{0,C}(m))\prod_{\{j\in J; S(j)\subset C\}} \Gamma(\alpha g_{j,C}(m))}{\Gamma(\alpha)\prod_{S\in \mathcal{S}} \left[\Gamma(\alpha g_{0,S}(m))\prod_{\{j\in J; S(j)\subset S\}} \Gamma(\alpha g_{j,S}(m))\right]^{\nu(S)}}.\label{formulai}
\end{eqnarray}
 and
  \begin{eqnarray}
  \lim_{\alpha\rightarrow}\alpha^{|J|}I(m,\alpha)&=&J_C(m)\nonumber\\
  &=&\frac{\prod_{S\in \mathcal{S}} \left[ g_{0,S}(m)\prod_{\{j\in J; S(j)\subset S\}}  g_{j,S}(m)\right]^{\nu(S)}}{\prod_{C\in \mathcal{C}} g_{0,C}(m)\prod_{\{j\in J; S(j)\subset C\}}  g_{j,C}(m)}\label{limit}
\end{eqnarray}  
  \end{theorem}
 
\begin{cor} 
\label{facedecomp}
In the case of a hierarchical model associated to a decomposable graph,  all the facets of $\overline{C}$ are of the type $F(j_0,C)$ described in Theorem \ref{generalfacets}, with $j_0\in J,$ with  $C$ in the set $\mathcal{C}$ of cliques and $S(j_0)\subset C.$
\end{cor}

 \begin{pff}
  We know from Theorem \ref{generalfacets} that the affine forms in the denominator of  $J_C(m)$ in (\ref{limit}) define facets of $\overline{C}$. From Theorem \ref{jmlambda1}, we know that they are the only ones.
 \end{pff}
In fact we conjecture, as mentioned in the introduction, that if a model is such that the only facets of $\overline{C}$ are of the type given in Theorem \ref{generalfacets}, then it is  a decomposable graphical model. 
 
 \vspace{4mm}\noindent {\it Example.} If $V=\stackrel{a}{\bullet}-\stackrel{b}{\bullet} -\stackrel{c}{\bullet}$ and if $I=\{0,1,2\}\times \{0,1\}\times \{0,1\}$ we have 
 \begin{eqnarray*}
 g_{0,bc}(m)&=&1-m_{001}-m_{010}+m_{011}\\ 
 g_{001,bc}(m)&=& m_{001}-m_{011}\\
 g_{010,bc}(m)&=&  m_{010}-m_{011}\\
 g_{011,bc}(m)&=& m_{011}\\
 g_{0,ab}(m)&=&1-m_{100}-m_{200}-m_{010}+m_{110}+m_{210}\\ 
 g_{100,ab}(m)&=&m_{100}-m_{110}\\
 g_{200,ab}(m)&=&m_{200}-m_{210}\\
 g_{010,ab}(m)&=&m_{010}-m_{110}-m_{210}\\
  g_{110,ab}(m)&=&m_{110}\\
 g_{210,ab}(m)&=&m_{210}\\
  g_{0,b}(m)&=&1-m_{010}\\
 g_{010,b}(m)&=&m_{010}
 \end{eqnarray*}
 In this case $I(m,\alpha)$ is a quotient: the numerator is  the product of  10  gamma functions and the denominator is $\Gamma(\alpha)\Gamma(\alpha(1-m_{010}))\Gamma(\alpha m_{010})$. As a consequence $J_C(m)$ is 
 $$\frac{g_{0,b}(m)g_{010,b}(m)}{g_{0,bc}(m)g_{001,bc}(m)g_{010,bc}(m)g_{011,bc}(m)g_{0,ab}(m)g_{100,ab}(m)g_{200,ab}(m)g_{010,ab}(m)g_{110,ab}(m)g_{210,ab}(m)}$$

\subsection{Facets of $\overline{C}$ when the model is binary and the model is governed by a cycle}
\label{ex:polytope}
For the sake of completion and for the convenience of the reader, we recall some known results giving 
 the facets of the polytope $\overline{C}$ when the model is hierarchical, binary and governed by a cycle $G$ of order $n\geq 3$.  The reader is referred to \cite{dl95} and \cite{hsul02} and some references within for an explicit description of these facets. In this subsection, we will simply translate the equation of the facets given in these papers in our own coordinates. The results are given in the following theorem. The coordinates of $m\in \mbox{R}^J$ will be denoted $m_v$ if they are indexed by a vertex $v\in V$ and by $m_e$ if they are indexed by an edge $e\in E$.
 \begin{theorem} \label{facetscycle} Let $G=(V,E)$ be a cycle of order $n\geq 3.$ Assume the hierarchical model is binary and governed by $G$, that is  ${\cal D}=\{v\in V,\;e\in E\}$. Then  the polytope $\overline{C}$ is defined by the following equations and the facets are defined by the corresponding equalities:
 \begin{enumerate}
  \item for any edge $(a,b)\in E$,
  \begin{eqnarray}
  &&m_{ab}\geq 0\label{c1}\\
   &&m_a-m_{ab}\geq 0\label{c2}\\
    &&m_b-m_{ab}\geq 0\label{c3}\\
     &&1-m_a-m_b+m_{ab}\geq 0\;,\label{c4}\\
  \end{eqnarray}
 \item for any subset $F\subseteq E$ with odd cardinality $|F|$,
 \begin{eqnarray}
 \hspace{1cm}\sum_{(a,b)\in F}(m_a+m_b-2m_{ab})-\Big(\sum_{v\in V}m_v-\sum_{e\in E}m_e\Big)\leq \frac{|F|-1}{2}\;.\label{c5}
 \end{eqnarray}
 \end{enumerate}
 The total number of facets for the polytope $\overline{C}$ of the model governed by the cycle of order $n$ is $F_n=\sum_{k\in N,k\;\mbox{odd},k\leq n}\left(\begin{array}{c} n\\ k\end{array}\right)$.
 \end{theorem}
 We see that the facets given by the first four equations are those described in Theorem \ref{generalfacets} corresponding to the cliques $\{(a,b)\in E\}$ while the others are specific to models governed by a cycle. We illustrate this theorem in the case of the cycles of order 3,4 and 5.
 We will not repeat the facets (\ref{c1})-(\ref{c4}) common to all hierarchical models. We will give the facets of type (\ref{c5}) only.
 
\noindent For $n=3$, let $V=\{a,b,c\}$ and $E=\{(a,b),(b,c),(c,a)\}$, the four facets of type (\ref{c5}) are
\begin{eqnarray}
&&1-m_a-m_b-m_c+m_{ab}+m_{bc}+m_{ac}\geq 0\nonumber\\
&&m_{ab}+m_c-m_{bc}-m_{ac}\geq 0\label{c3-1}
\end{eqnarray}
and the other two facets obtained from (\ref{c3-1}) by permutations of the edges of $G$.
\vspace{2mm}

\noindent For $n=4$, let $V=\{a,b,c,d\}$ and $E=\{(a,b),(b,c),(c,d), (d,a)\}$, the eight facets of type (\ref{c5}) are

\begin{eqnarray}
&&1-m_a-m_b-m_c-m_d+m_{ab}+m_{bc}+m_{cd}+m_{da}\geq 0\nonumber\\
&&m_c+m_d+m_{ab}-m_{bc}-m_{cd}-m_{da}\geq 0\label{c4-1}
\end{eqnarray}
and the other three facets obtained from (\ref{c4-1}) by permutations of the edges of $G$.
\vspace{2mm}

\noindent For $n=5$, let $V=\{a,b,c,d,e\}$ and $E=\{(a,b),(b,c),(c,d), (d,e), (e,a)\}$, sixteen facets of type (\ref{c5}) are
\begin{eqnarray}
&&m_{ab}+m_c+m_d+m_e-m_{bc}-m_{cd}-m_{de}-m_{da}\geq 0\label{c5-1}\\
&&1-m_a-m_b+m_{ea}+m_{ab}+m_{bc}+m_d-m_{cd}-m_{ed}\geq 0\label{c5-2}\\
&&1-m_d+m_{ab}+m_{cd}+m_{de}-m_{bc}-m_{ae}\geq 0\label{c5-3}\\
&&2-m_a-m_b-m_c-m_d-m_e+m_{ab}+m_{bc}+m_{cd}+m_{de}+m_{ea}\geq 0\nonumber
\end{eqnarray}
and the other three facets obtained from each of (\ref{c5-1}),(\ref{c5-2}) and (\ref{c5-3})  by permutations of the edges of $G$.

\section{Conclusion} 
The main contribution of this paper is the description of the behaviour of the Bayes factor as $\alpha\rightarrow 0$. We have shown that, in this study, the important concept is the dimension of the face to which the data belongs rather than the dimension of the model. We have identified the role of the open convex polytope $C$ and the function $J_C(.)$ It is not surprising that $\overline{C}$, the convex hull of the support of the generating measure of the multinomial for the hierarchical model, plays an important role. The multinomial for a loglinear hierarchical model is a natural exponential family and the role of $C$ which is the domain of the mean is well-known. The  set $C$ is also of prime importance  in the study of the existence of the maximum likelihood estimate of the parameter (see for example Eriksson et al. \cite{efrs06} or Geiger et al. \cite{gms06} or  Rinaldo \cite{rinaldo}).  However, the role of  the characteristic function $J_C(.)$ of $C$ has only been uncovered now in the study of the Bayes factor and we can add $J_C$ to the toolkit of exponential families. We note that all the limit theorems in Section 3 are valid for any natural exponential family such that the convex hull of the support of its generating measure is a bounded convex polytope but are not be immediately applicable to a family of distributions such as the Poisson where $C$ is not bounded. This is the topic of further work.

 A secondary contribution of this paper is our results on the identification of the facets of a polytope. We have two new results for polychotomous models (i.e. not necessarily binary): the first giving a particular category of facets common to all hierarchical models, the second giving the complete set of facets for decomposable models. 
 
 We have also extended the results of \cite{sj02} to the case of any two decomposable models, thus allowing the practitioner to predict the behaviour of the Bayes factor without using the concept of face or facets of a polytope.

 \section{Acknowlegments} We are grateful to Steffen Lauritzen who asked the question which motivated this paper, that is, what happens to the Bayes factor when $\alpha \rightarrow 0$, and who told us about \cite{sj02}. We thank Seth Sullivant
for pointing out the work of \cite{dl95} and \cite{hsul02},  Jean-Baptiste Hiriart-Urruty for  refering us to \cite{hiriart}, Alexander Barvinok who directed us to the second formula of (\ref{eq:J}), Mathieu Meyer who suggested a proof of Theorem 3.2 much shorter than our initial one, Monique Laurent who helped us with Theorem 5.3 and an anonymous contributor ('Meu') of the forum Les-Mathematiques.net who gave us an outline of the proof of Lemma \ref{meu}.

\end{document}